\documentclass[a4wide, 12pt]{article}
\usepackage{amsmath, amsfonts, amssymb, amsthm}
\usepackage{indentfirst}

\usepackage{fancyheadings}

\newtheorem{theorem}{Theorem}
\newtheorem{lemma}{Lemma}
\newtheorem{statement}{Statement}

\newenvironment{Proof}
{\par\noindent{\bf Proof.}}
{\hfill$\scriptstyle\square$}

\newcommand{\Z}{\mathbb{Z}}
\newcommand{\N}{\mathbb{N}}
\newcommand{\R}{\mathbb{R}}
\newcommand{\Pl}{\mathcal{P}}
\newcommand{\A}{\mathbb{A}}

\title{On points with algebraically conjugate coordinates close to smooth curves}
\author{V. Bernik, F. G\"otze, A. Gusakova}

\begin{document}

\maketitle

\begin{abstract}

Let $y=f(x)$ be a continuous differentiable function on an interval $J\subset\R$. In this paper we show that for any $n\in\N$, $n\ge 2$, sufficiently large integer $Q$ and a real $0<\lambda<\frac34$ there exists a positive value $c(n,f,J)$ such that all strips $L_{J}(Q,\lambda)=\left\{(x_1,x_2)\in\R^2:\left|x_2-f(x_1)\right|\ll Q^{-\lambda}, x_1\in J\right\}$ contain at least $c(n, f, J)Q^{n+1-\lambda}$ points $\boldsymbol{\alpha}=(\alpha_1,\alpha_2)$ with algebraically conjugate coordinates which minimal polynomial $P$ satisfies $\deg P\leq n$, $H(P)\leq Q$. The proof is based on a metric theorem on the measure of the set of vectors $(x_1,x_2)$ lying in a rectangle $\Pi$ of dimensions $\asymp Q^{-s_1}\times Q^{-s_2}$ with $|P(x_1)|, |P(x_2)|$ bounded from above and $|P'(x_1)|,|P'(x_2)|$ bounded from below, where $P$ is a polynomial of degree $\deg P\leq n$ and height $H(P)\leq Q$. This theorem is a generalization of a result obtained by V. Bernik, F. G\"otze and O. Kukso for $s_1=s_2=\frac12$ and $\lambda = \frac12$ \cite{BernikGoetzeKukso14} .


\end{abstract}

\section{Introduction}

Let $Q$ be a sufficiently large number. We denote by $\Pl_n(Q)$ the following class of polynomials:
\[
\Pl_n(Q)=\{P\in\Z[t]:\deg P\leq n, H(P)\leq Q\},
\]
where $H(P)=\max\limits_{0\leq j \leq n}{|a_j|}$ denotes the height of an integer polynomial $P(t)=a_nt^n+\ldots +a_1t+a_0$. 

The point $\boldsymbol{\alpha}=(\alpha_1,\alpha_2)$ is called an {\it algebraic point} if $\alpha_1$ and $\alpha_2$ are roots of the same polynomial $P\in\Z[t]$. The polynomial $P$ of smallest degree such that $P(\alpha_1)=P(\alpha_2)=0$ and $\gcd\left(|a_n|,\ldots,|a_0|\right)=1$ is called the minimal polynomial of the algebraic point $\boldsymbol{\alpha}$. Denote by $\deg(\boldsymbol{\alpha})=\deg P$ the degree of the algebraic point $\boldsymbol{\alpha}$, and by $H(\boldsymbol{\alpha})=H(P)$ the height of the algebraic point $\boldsymbol{\alpha}$. Define the following sets: $\A_n^2(Q)$ is the set of algebraic points $\boldsymbol{\alpha}$ of degree at most $n$ and of height at most $Q$; $\A_n^2(Q, D)=\A_n^2(Q)\cap D$ is the set of algebraic points $\boldsymbol{\alpha}\in\A_n^2(Q)$ lying in a domain $D\subset \R^2$. Denote by $\#S$ the cardinality of a finite set $S$, by $\mu_1 S$ the Lebesgue measure of a measurable set $S\subset \R$ and by $\mu_2 S$ the Lebesgue measure of a measurable set $S\subset \R^2$. Further, denote by $c_j> 0$, $j\in\N$, positive values which do not depend on $H(P)$ or $Q$. We are also going to use the Vinogradov symbol $A\ll B$, which means that there exists a value $c>0$ such that $A\leq c\cdot B$ and $c$ doesn't depend on $B$.

An important and interesting topic in the theory of Diophantine approximation is the distribution of algebraic numbers \cite{Sprindzuk67, Cassels57, Schmidt80, BeresnevichBernikGoetze10}. In this paper we consider problems related to the distribution of algebraic points in domains of small measure and the distribution of algebraic points near smooth curves. 

Consider rectangles $\Pi = I_1\times I_2$ where $\mu_1 I_1=c_{1,1}\cdot Q^{-s_1}$ and $\mu_1 I_2=c_{1,2}\cdot Q^{-s_2}$ under the conditions $0< s_1+s_2\leq 1$, $s_1, s_2 < 1$, $\Pi\cap\left\{(x_1,x_2)\in\R^2:\ |x_1-x_2|<\varepsilon\right\}=\emptyset$ and $c_{1,1}c_{1,2} \ge c_0$.
The condition $|x_1-x_2|>\varepsilon$ means that we exclude from consideration a strip $F$ of small measure such that the coordinates $(x_1,x_2)\in F$ are well approximated by points of form $(\alpha,\alpha)$.

We can prove the following theorem.

\begin{theorem}\label{th1}
For any rectangle $\Pi=I_1\times I_2$ satisfying the following conditions:

1. $\mu_1 I_i = c_{1,i} Q^{-s_i}$ where $s_i < 1$ and $0< s_1+s_2\leq 1$, $i=1,2$;

2. $\Pi\cap\left\{(x_1,x_2)\in\R^2:\ |x_1-x_2|<\varepsilon\right\}=\emptyset$;

3. $c_{1,1}c_{1,2}> c_0(n,\varepsilon,\mathbf{d})>0$  for $s_1+s_2=1$, where $\mathbf{d}=(d_1,d_2)$ is the midpoint of $\Pi$;

\noindent there exists a constant $c_2=c_2(n,\varepsilon,\mathbf{d})>0$, such that
\[
\#\A_n^2(Q,\Pi)\ge c_2 Q^{n+1}\mu_2\Pi, 
\]
for $Q>Q_0(n,\varepsilon,\mathbf{d},\mathbf{s})$.
\end{theorem}

For $s_1+s_2>1$, we can find a rectangle $\Pi$ such that the statement of Theorem 1 does not hold. The example of such rectangle is $\Pi=\left(0, 0.5Q^{-1}\right)\times \left(0, 0.5\right)$. It is easy to prove \cite{BernikGoetze14} that the interval $\left(0, 0.5Q^{-1}\right)$ doesn't contain algebraic numbers of any degree and height $\leq Q$. Let us introduce some restrictions on the domains to be used in the following proofs. 

Consider a square $\overline{\Pi}=I_1\times I_2$ of size $\mu_1 I_1 = \mu_1 I_2 = c_3Q^{-s}$ such that $\frac12< s<\frac34$. Given positive $u_1,u_2$ under the condition $u_1+u_2=1$ let us say that the square $\overline{\Pi}$ is {\it $(u_1,u_2)$- ordinary} square if it doesn't contain points $(x_1',x_2')\in \R^2$ such that there exists a polynomial $P\in\Pl_2(Q)$ of the form $P(t)=b_2t^2+b_1t+b_0$ satisfying the system of inequalities
\begin{equation}\label{eq1}
\begin{cases}
|P(x_i')|\ll Q^{-u_i},\quad i=1,2,\\
|b_2|<Q^{s-\frac12}.
\end{cases}
\end{equation}
Otherwise, the square $\overline{\Pi}$ is going to be called {\it $(u_1,u_2)$- special}.

For {\it $(u_1,u_2)$-ordinary} squares, the following result holds.

\begin{theorem}\label{th2}
For any {\it $\left(\frac12,\frac12\right)$-ordinary} square $\overline{\Pi}=I_1\times I_2$ under the following conditions:

1. $\mu_1 I_i = c_3 Q^{-s}$, where $\frac12< s< \frac34$;

2. $\overline{\Pi}\cap\left\{(x_1,x_2)\in\R^2:\ |x_1-x_2|<\varepsilon\right\}=\emptyset$;

3. $c_3>c_0(n,\varepsilon,\mathbf{d})>0$, where $\mathbf{d}=(d_1,d_2)$ is the midpoint of $\overline{\Pi}$;

\noindent there exists a constant $c_4=c_4(n,\varepsilon, \mathbf{d})>0$, such that
\[
\#\A_n^2(Q,\overline{\Pi})\ge c_4 Q^{n+1}\mu_2\overline{\Pi}
\]
for $Q>Q_0(n,\varepsilon,\mathbf{d},s)$.
\end{theorem}

Another interesting and important topic is the distribution of algebraic points near smooth curves. The result presented in this paper is a natural generalization of problems related to distribution of rational points near smooth curves \cite{Huxley, BereDickVel, BeresnevichBernikGoetze10, Bugeaud02, BugeaudMignotte04, Evertse04}. In 2014 a lower bound for the number of algebraic points lying at a distance of at most $Q^{-\lambda}$, $0<\lambda<\frac12$, from a smooth curve was obtained by V. Bernik, F. G\"otze and O. Kukso \cite{BernikGoetzeKukso14}. We improve on this result and obtain an identical estimate for $0<\lambda<\frac34$.

\begin{theorem}\label{th3}
Let $y=f(x)$ be a continuous differentiable function on an interval $J=[a,b]$ such that $\sup\limits_{x\in J}{|f'(x)|}:=c_5 <\infty$ and $\#\left\{x\in\R: f(x)=x\right\}<\infty$. Denote by $L_{J}(Q,\lambda)$ the following set:
\[
L_{J}(Q,\lambda)=\left\{(x_1,x_2)\in\R^2:|x_2-f(x_1)|<\left(\textstyle\frac12+c_5\right)\cdot c_3Q^{-\lambda},\quad x_1\in J\right\},
\]
for $0<\lambda<\frac34$. Then there exists a positive value $c_6(J,f,n)>0$ such that
\[
\#\left\{\A_n^2(Q)\cap L_J\left(Q,\lambda\right)\right\}\ge c_6Q^{n+1-\lambda}
\]
for $Q>Q_0(J,f,n,\lambda)$.
\end{theorem}

\section{Auxiliary statements}

For a polynomial $P$ with roots $\alpha_1,\alpha_2,\ldots,\alpha_n$, let 
\[
S(\alpha_i) = \left\{x\in\R : |x-\alpha_i| = \min\limits_{1\le j\le n} |x-\alpha_j|\right\}.
\]
From now on, we assume that the roots of the polynomial $P$ are sorted by distance from $\alpha_i=\alpha_{i,1}$:
\[
|\alpha_{i,1}-\alpha_{i,2}|\le|\alpha_{i,1}-\alpha_{i,3}|\le\ldots\le|\alpha_{i,1}-\alpha_{i,n}|.
\]

\begin{lemma}\label{lm1}
Let $x\in S(\alpha_i)$. Then
\begin{align}
&|x-\alpha_i| \le n\cdot \frac{|P(x)|}{|P'(x)|},\quad |x-\alpha_i| \le 2^{n-1}\cdot \frac{|P(x)|}{|P'(\alpha_i)|},\\
&|x-\alpha_i| \le \min\limits_{1\le j\le n} \left(2^{n-j} \frac{|P(x)|}{|P'(\alpha_i)|}
|\alpha_{i}-\alpha_{i,2}|\ldots|\alpha_{i}-\alpha_{i,j}| \right)^{1/j}.
\end{align}
\end{lemma}

The first inequality follows from the identity
\[
|P'(x)||P(x)|^{-1} = \sum\limits_{j=1}^n |x-\alpha_j|^{-1}.
\]
For a proof of the second and the third inequalities see \cite{Sprindzuk67}, \cite{Bernik83}.

\begin{lemma}\label{lm2}
Let $I$ be an interval, and let $A\subset\R$ be a measurable set,
$A\subset I$, $\mu_1 A\ge \frac12 \mu_1 I$. If for some $v > 0$ and all $x\in A$ the inequality
$|P(x)|<c_7Q^{-v}$, where $v>0$, holds, then
\[
|P(x)| < 6^n(n+1)^{n+1}c_7 Q^{-v}
\]
for all points $x\in I$, where $n=\deg P$.
\end{lemma}

The proof of this lemma can be found in \cite{Bernik80}.

\begin{lemma}\label{lm3}
Let $\delta$, $\eta_1$, $\eta_2$ be real positive numbers,
and let $P_1,P_2 \in \Z[t]$ be a co-prime polynomials of degrees at most $n$ such that
\[
\max\left(H(P_1), H(P_2)\right) < K,
\]
where $K>K_0(\delta)$. Let
$J_1,J_2\subset \R$ be intervals of sizes $\mu J_1=K^{-\eta_1}$, $\mu J_2=K^{-\eta_2}$. If for some $\tau_1, \tau_2>0$ and for all $(x_1,x_2)\in J_1\times J_2$, the inequalities
\[
\max\left(|P_1(x_i)|, |P_2(x_i)|\right) < K^{-\tau_i},\quad i=1,2,
\]
hold, then
\begin{equation}\label{eq0}
\tau_1+\tau_2+2 + 2\max(\tau_1+1-\eta_1, 0) + 2\max(\tau_2+1-\eta_2, 0) < 2n+\delta.
\end{equation}
\end{lemma}

The proof of this lemma can be found in \cite{Pereverzeva}.

\begin{lemma}\label{lm4}
Let $P\in\Z[t]$ be a reducible polynomial, $P=P_1\cdot P_2$,
$\deg P = n\ge 2$. Then there exist $c_8,c_{9}>0$ such that
\[
c_{8}H(P)< H(P_1) H(P_2) <c_{9}H(P).
\]
\end{lemma}

The proof of Lemma \ref{lm4} can be found, for example, in 
\cite{Sprindzuk67}.

\section{Proof of Theorem 1}

Before we start it should be noted that there exists a constant $h_n=h_n(\mathbf{d})>0$ such that for every point $(x_1,x_2)\in\Pi$ and every $\mathbf{v}=(v_1,v_2)$ with $v_1+v_2=n-1$ there exists a polynomial $P\in\Pl_n(Q)$ satisfying the inequalities:
\[
|P(x_i)|< h_n\cdot Q^{-v_i},\quad i=1,2,
\]
for $Q>Q_0$. This simple fact follows from Dirichlet's principle and estimates $\# \Pl_n(Q)>2^nQ^{n+1}$ and $|P(x_i)|< \left(\left(|d_i|+1\right)^{n+1}-1\right)|d_i|^{-1}\cdot Q$, where $\mathbf{d}=(d_1,d_2)$ is the midpoint of $\Pi$.

To prove Theorem \ref{th1}, we are going to rely on the following Lemma \ref{lm5}. 
 
\begin{lemma}\label{lm5}
For all rectangles $\Pi=I_1\times I_2$ under the conditions:

1. $\mu_1 I_i = c_{1,i} Q^{-s_i}$ where $s_i < 1$ and $0< s_1+s_2\leq 1$, $i=1,2$;

2. $\Pi\cap\left\{(x_1,x_2)\in\R^2:\ |x_1-x_2|<\varepsilon\right\}=\emptyset$;

3. $c_{1,1}c_{1,2}> c_0(n,\varepsilon,\mathbf{d})>0$  for $s_1+s_2=1$, where $\mathbf{d}=(d_1,d_2)$ is the midpoint of $\Pi$;

\noindent let $L=L(Q,\delta_n, \mathbf{v},\Pi)$ be the set of points $(x_1,x_2)\in\Pi$ such that there exists a polynomial $P\in\Pl_n(Q)$ satisfying the following system of inequalities:
\begin{equation}\label{eq5}
\begin{cases}
|P(x_i)|< h_n\cdot Q^{-v_i},\quad v_i>0,\\
\min\limits_i\left\{|P'(x_i)|\right\}<\delta_n\cdot Q,\\
v_1+v_2=n-1,\quad i=1,2.
\end{cases}
\end{equation}
Then for $\delta_n\leq \delta_0(n,\varepsilon,\mathbf{d})$ and $Q>Q_0(n,\varepsilon,\mathbf{s},\mathbf{v},\mathbf{d})$, the estimate
\[
\mu_2 L<\textstyle\frac14 \mu_2\Pi
\]
holds. 
\end{lemma}

\begin{Proof}
Denote by $L_{1}$ the set of points $(x_1,x_2)\in\Pi$ such that the system of inequalities (\ref{eq5}) has a solution in irreducible polynomials $P\in\Pl_n(Q)$ under condition $|P'(x_1)|<\delta_n\cdot Q$, by $L_2$ the set of points $(x_1,x_2)\in\Pi$ such that the system of inequalities (\ref{eq5}) has a solution in irreducible polynomials $P\in\Pl_n(Q)$ under condition $|P'(x_2)|<\delta_n\cdot Q$ and by $L_3$ the set of points $(x_1,x_2)\in\Pi$ such that the system of inequalities (\ref{eq5}) has a solution in reducible polynomials $P\in\Pl_n(Q)$. Thus, $L=L_1\cup L_2\cup L_3$. 

Let us estimate the measure of $L_1$. The main idea is to split the range of the possible values of $|P'(x_i)|$, $|P'(\alpha_i)|$, where $x_i\in S(\alpha_i)$, $i=1,2$ into a total of $r=r(n)=(n-1)^2$ sub-ranges and consider them separately.

Without loss of generality, we will assume that $|d_1|< |d_2|$. Let us show that the inequality
\begin{equation}\label{eq6}
|P'(x_i)|\ge 2c_{10}\cdot Q^{\frac12-\frac{v_i}{2}}
\end{equation}
yields the following bounds on $P'(\alpha_i)$:
\[
\textstyle\frac12|P'(x_i)| \leq |P'(\alpha_i)|\leq 2|P'(x_i)|,
\]
where $c_{10}=n(n-1)\cdot\max\{h_n,1\}\cdot\left(3\max\left\{1,|d_2|\right\}\right)^{n-1}\cdot \left(1+|d_2|^{-1}\right)$. Let us write a Taylor expansion of $P'(t)$:
\begin{equation}\label{eq7}
P'(x_i)=P'(\alpha_i)+\textstyle\frac12P''(\alpha_i)(x_i-\alpha_i)+\ldots+\textstyle\frac{1}{(n-1)!}P^{(n)}(\alpha_i)(x_i-\alpha_i)^{n-1}.
\end{equation}
Using Lemma \ref{lm1} and the estimates (\ref{eq5}) for $Q>Q_0$, we have:
\[
|x_i-\alpha_i|\leq nh_nc_{10}^{-1} \cdot Q^{-\frac{v_i+1}{2}}< \max\left\{1,|d_2|\right\}\cdot Q^{-\frac{v_i+1}{2}}.
\]
Then, for $s_i>0$ and $Q>Q_0$ we get $|x_i-d_1|<1/2$ and thus:
\[
|\alpha_i|\leq |x_i|+\textstyle\frac12< |d_2|+1.
\]
From this estimates we obtain the following inequality for every term in (\ref{eq7}):
\begin{multline*}
\left|\textstyle\frac{1}{(k-1)!}P^{(k)}(\alpha_i)(x_i-\alpha_i)^{k-1}\right|< C^{k-1}_{n-1}\cdot \textstyle\frac{n(n+1-k)(|d_2| + 1)^{n-k+1}}{|d_2|}\cdot \max\left\{1,|d_2|\right\}^{k-1}\cdot Q^{1-\frac{(k-1)(1+v_i)}{2}}\leq\\
\leq C^{k-1}_{n-1}\cdot \textstyle\frac{n(n-1)(|d_2| + 1)^{n-k+1}}{|d_2|}\cdot \max\left\{1,|d_2|\right\}^{k-1}Q^{\frac12-\frac{v_i}{2}},
\end{multline*}
for $k\ge 2$. Thus, the estimate
\begin{multline*}
\left|\textstyle\frac12P''(\alpha_i)(x_i-\alpha_i)+\ldots+\textstyle\frac{1}{(n-1)!}P^{(n)}(\alpha_i)(x_i-\alpha_i)^{n-1}\right|< \\
< n(n-1)\left(3\max\left\{1,|d_2|\right\}\right)^{n-1}\cdot \left(1+|d_2|^{-1}\right)\cdot Q^{\frac12-\frac{v_i}{2}}<\textstyle\frac12\cdot |P'(x_i)|
\end{multline*}
holds. By substituting these inequality to (\ref{eq7}) we get 
\[
\textstyle\frac12\cdot |P'(x_i)|\leq|P'(\alpha_i)|\leq 2|P'(x_i)|.
\]

This means that $|P'(\alpha_i)|\in T_i$, where 
\[
T_1=\left[c_{10}\cdot Q^{\frac12-\frac{v_1}{2}}; 2\delta_n\cdot Q\right), \qquad\qquad T_2=\left[c_{10}\cdot Q^{\frac12-\frac{v_2}{2}}; n\cdot\textstyle\frac{\left(|d_2|+1\right)^{n}-1}{|d_2|}\cdot Q\right)
\]
if the inequalities (\ref{eq6}) hold. Let us divide the intervals $T_i$ into sub-intervals $T_{i,j}=\left[d_{j,i}Q^{t_{j,i}};d_{j-1,i}Q^{t_{j-1,i}}\right)$, $2\leq j\leq n$, where 
\[
t_{k,i}=
\begin{cases}
1,\quad k=1,\\
\frac12-\frac{(k-1)v_i}{2(n-1)},\quad 2\leq k\leq n,
\end{cases}\qquad
d_{k,i}=
\begin{cases}
2\delta_n,\quad k=1,i=1,\\
n\cdot\textstyle\frac{\left(|d_2|+1\right)^{n}-1}{|d_2|},\quad k=1,i=2,\\
1,\quad 2\leq k\leq n-1,\\
c_{10},\quad k=n,
\end{cases}
\]
Now we are going to consider the following cases:
\begin{itemize}
\item the case of polynomials of the second degree $n=2$ (see Section \ref{sec_1});
\item the case of irreducible polynomials:
\subitem $|P'(\alpha_1)|\in T_{1,j_1}$, $|P'(\alpha_2)|\in T_{2,j_2}$, where $1\leq j_1,j_2\leq n-1$ (see Section \ref{sec_2});
\subitem $|P'(\alpha_1)|\in T_{1,n}$, $|P'(\alpha_2)|\in T_{2,n}$ (see Section \ref{sec_3});
\subitem $|P'(x_1)|\leq 2c_{10}Q^{\frac12-\frac{v_1}{2}}$, $|P'(x_2)|\leq 2c_{10}Q^{\frac12-\frac{v_2}{2}}$ (see Section \ref{sec_4});
\subitem $|P'(\alpha_1)|\in T_{1,j_1}$, $|P'(\alpha_2)|\in T_{2,n}$ or $|P'(\alpha_1)|\in T_{1,n}$, $|P'(\alpha_2)|\in T_{2,j_2}$, where $1\leq j_1,j_2\leq n-1$ (see Section \ref{sec_5});
\subitem $|P'(\alpha_1)|\in T_{1,j_1}$, $|P'(x_2)|\leq 2c_{10}Q^{\frac12-\frac{v_2}{2}}$ or $|P'(x_1)|\leq 2c_{10}Q^{\frac12-\frac{v_1}{2}}$, $|P'(\alpha_2)|\in T_{2,j_2}$, where $1\leq j_1,j_2\leq n$ (see Section \ref{sec_5});
\item the case of reducible polynomials (see Section \ref{sec_6}).
\end{itemize}

We are going to use induction on the degree $n$. Let us prove the following statement, which will serve as the base of induction.

\subsection{The base of induction: polynomials of the second degree.}\label{sec_1}

\begin{statement}\label{st1}
For all rectangles $\Pi$ under the conditions 1--- 3 let $L_{2,2}=L_{2,2}(Q,\delta_2,\boldsymbol{\gamma}_2,\Pi)$ be the set of points $(x_1,x_2)\in\Pi$ such that there exists a polynomial $P\in\Pl_2(Q)$ satisfying the system of inequalities
\begin{equation}\label{eq9}
\begin{cases}
|P(x_i)|< h_2\cdot Q^{-\gamma_{2,i}},\quad \gamma_{2,i} > 0,\\
\min\limits_i\left\{|P'(x_i)|\right\}<\delta_2\cdot Q,\\
\gamma_{2,1}+\gamma_{2,2}=1,\quad i=1,2.
\end{cases}
\end{equation}
Then for any $r>0$ and for $\delta_2<\delta_0(r,\varepsilon,\mathbf{d})$ and $Q>Q_0(n,\varepsilon,\mathbf{s},\boldsymbol{\gamma}_2,\mathbf{d})$, the estimate
\[
\mu_2 L_{2,2}<\textstyle\frac{1}{4r}\cdot\mu_2\Pi
\]
holds.
\end{statement}

\begin{Proof}
Let $P(t)$ be a polynomial of the form $b_2t^2+b_1t+b_0$. Let us estimate the values $|P'(\alpha_1)|$ and $|P'(\alpha_2)|$. By the third inequality of Lemma \ref{lm1}, for every polynomial $P$ satisfying the inequalities (\ref{eq9}) at a point $(x_1,x_2)\in\Pi$, we have the following estimates:
\begin{equation}\label{eq10}
|x_i-\alpha_i|<\left(|P(x_i)||b_2|^{-1}\right)^{1/2} < h_2^{1/2} Q^{-\frac{\gamma_{2,i}}{2}} < \textstyle\frac{\varepsilon}{8},
\end{equation}
for $Q>Q_0$ and $x_i\in S(\alpha_i)$, $i=1,2$.

From (\ref{eq10}) and condition 2 we obtain that
\[
|\alpha_1-\alpha_2|>|x_1-x_2|-|x_1-\alpha_1|-|x_2-\alpha_2|>\textstyle\frac{3}{4}\cdot\varepsilon
\]
and
\[
|\alpha_1-\alpha_2|<|x_1|+|x_2|+|x_1-\alpha_1|+|x_2-\alpha_2|< |d_1|+|d_2|+1 + \textstyle\frac{\varepsilon}{4}.
\]
This leads to the following lower bounds for $|P'(\alpha_i)|$:
\begin{equation}\label{eq0011}
\left(|d_1|+|d_2|+1+\textstyle\frac{\varepsilon}{4}\right)\cdot|b_2|>|P'(\alpha_i)|=\sqrt{D}=|b_2|\cdot|\alpha_1-\alpha_2|>\textstyle\frac{3}{4}\cdot\varepsilon\cdot|b_2|,
\end{equation}
where $D$ is the discriminant of the polynomial $P$.
The inequalities (\ref{eq10}) also yield upper bounds for $|P'(x_i)|$:
\begin{equation}\label{eq011}
|P'(x_i)|\leq |b_2|\cdot\left(|\alpha_1-x_i|+|\alpha_2-x_i|\right)\leq \left(|d_2|+1+\textstyle\frac{\varepsilon}{4}\right)\cdot|b_2|.
\end{equation}

Now upper bounds for $|P'(\alpha_i)|$ can be obtained from the Taylor expansion of the polynomial $P'$:
\begin{equation}\label{eq11}
|P'(\alpha_i)|\leq|P'(x_i)|+|P''(x_i)|\cdot|x_i-\alpha_i|\leq |P'(x_i)|+\textstyle\frac{\varepsilon}{2}\cdot|b_2|.
\end{equation}

Then, the estimates (\ref{eq0011}), (\ref{eq11}) mean that
\begin{equation}\label{eq12}
|b_2|<4\varepsilon^{-1}\cdot\min\limits_i\left\{|P'(x_i)|\right\}<4\delta_2\varepsilon^{-1}Q.
\end{equation}

From Lemma \ref{lm1} and the estimates (\ref{eq0011}) it follows that the set $L_{2,2}$ is contained in a union  
$\bigcup\limits_{P\in\Pl_2(Q)}{\sigma_{P}}$, where
\[
\sigma_{P}=
\left\{(x_1,x_2)\in\Pi:\quad
|x_i-\alpha_i|<2h_2\varepsilon^{-1}Q^{-\gamma_{2,i}}|b_2|^{-1},i=1,2\right\}.
\]
Simple calculations show that the measure of the set $\sigma_{P}$ is lower than the measure of the rectangle $\Pi$:
\[
\mu_2\sigma_{P}\leq 2^4h_2^2\varepsilon^{-2}Q^{-1}|b_2|^{-2}< c_{1,1}c_{1,2}Q^{-1}=\mu_2\Pi
\]
for $c_{1,1}c_{1,2}>2^4h_2^2\varepsilon^{-2}$.

Let us estimate the measure of $L_{2,2}$:
\[
\mu_2 L_{2,2}\leq\mu_2\bigcup\limits_{P\in\Pl_2(Q)}{\sigma_{P}}\leq\sum\limits_{P\in\Pl_2(Q)}{\mu_2\sigma_{P}}\leq 2^4h_2^2\varepsilon^{-2}Q^{-1}\sum\limits_{\substack{b_2,b_1,b_0\leq Q: \\ P(t)=b_2t^2+b_1t+b_0,\\ \sigma_{P}\neq\emptyset}}{|b_2|^{-2}}.
\]
To do this, we need to estimate the number of polynomials $P\in\Pl_2(Q)$ such that the system (\ref{eq9}) holds for some point $(x_1,x_2)\in\Pi$, where $b_2$ is fixed. 

Let the inequalities (\ref{eq9}) hold for polynomial $P$ and point $(x_{0,1}, x_{0,2})\in\Pi$. Let us estimate the value of the polynomial $P$ at $d_i$. 
From the Taylor expansion of $P$, we have
\[
P(d_i)=P(x_{0,i}) +P'(x_{0,i})(x_{0,i}-d_i) +\textstyle\frac12 P''(x_{0,i})(x_{0,i}-d_i)^2.
\]
It means that $|P(d_i)|\leq |P(x_{0,i})|+|P'(x_{0,i})|\mu_1 I_i+|b_2|\left(\mu_1 I_i\right)^2$.
Thus, from (\ref{eq011}) for $Q>Q_0$ we can obtain the estimate 
\[
|P(d_i)|< |P(x_{0,i})|+c_{11}\cdot|b_2|\mu_1 I_i\leq 2c_{11}\cdot\max\{1,|b_2|\mu_1 I_i\}.
\]
Without loss of generality, let us assume that $\mu_1 I_1\leq \mu_1 I_2$.

Consider the system of equations
\begin{equation}\label{eq14}
\begin{cases}
b_2d_1^2+b_1d_1+b_0=l_1,\\
b_2d_2^2+b_1d_2+b_0=l_2
\end{cases}
\end{equation}
in three variables $b_2,b_1,b_0\in\Z$, where $|l_i|\leq 2c_{11}\cdot\max\{1,|b_2|\mu_1 I_i\}$, $i=1,2$.

Let us estimate the number of possible pairs $(b_1,b_0)$ such that the system (\ref{eq14}) is satisfied for a fixed $b_2$. To obtain this estimate, we consider the system of linear equations (\ref{eq14}) for two different combinations $b_2, b_{0,1},b_{0,0}$ and $b_2, b_{j,1},b_{j,0}$:
\[
\begin{cases}
b_2d_1^2+b_{0,1}d_1+b_{0,0}=l_{0,1},\\
b_2d_1^2+b_{j,1}d_1+b_{j,0}=l_{j,1},\\
b_2d_2^2+b_{0,1}d_2+b_{0,0}=l_{0,2},\\
b_2d_2^2+b_{j,1}d_2+b_{j,0}=l_{j,2}.
\end{cases}
\]
Subtracting the second equation from the first and the forth equation from the third leads to the following system in two variables $b_{0,1}-b_{j,1}$ and $b_{0,0}-b_{j,0}$:
\begin{equation}\label{eq15}
\begin{cases}
(b_{0,1}-b_{j,1})d_1+(b_{0,0}-b_{j,0})=l_{0,1}-l_{j,1},\\
(b_{0,1}-b_{j,1})d_2+(b_{0,0}-b_{j,0})=l_{0,2}-l_{j,2}.
\end{cases}
\end{equation}
The determinant of the system (\ref{eq15}) can be written as
\[
|\Delta|= 
\begin{vmatrix}
d_1 & 1\\
d_2 & 1
\end{vmatrix}
=|d_1-d_2|>\varepsilon>0.
\]
Since the determinant does not vanish, we can use Cramer's rule to solve the system (\ref{eq15}). Using the inequalities $|l_{0,i}-l_{j,i}|\leq 4c_{11}\cdot\max\{1,|b_2|\mu_1 I_i\}$, $i=1,2$, we estimate the determinant $\Delta_1$ as follows:
\[
|\Delta_1|\leq 8c_{11}\cdot\max\{1,|b_2|\mu_1 I_2\}.
\]

Hence by Cramer's rule we have 
\[
|b_{0,1}-b_{j,1}|\leq\frac{|\Delta_1|}{|\Delta|}\leq 8\varepsilon^{-1}c_{11}\cdot\max\{1,|b_2|\mu_1 I_2\}.
\]
This inequality means that all possible values of the coefficient $b_{1}$ lie in an interval $J_1$ of length $\mu_1 J_1=2^4\varepsilon^{-1}c_{11}\cdot\max\{1,|b_2|\mu_1 I_2\}$ centered at $b_{0,1}$. Since the values of the coefficient $b_1$ are integers, the number of these values does not exceed the measure of the interval $J_1$.

In addition, let us fix the value of the coefficient $b_1$. Choose a value  $b_1\in J_1$ and consider two different combinations $(b_2, b_1, b_{0,0})$ and $(b_2, b_1, b_{j,0})$. In this case, the system (\ref{eq14}) can be transformed as follows:
\[
\begin{cases}
|b_{0,0}-b_{j,0}|\leq 4c_{11}\cdot\max\{1,|b_2|\mu_1 I_1\},\\
|b_{0,0}-b_{j,0}|\leq 4c_{11}\cdot\max\{1,|b_2|\mu_1 I_2\}.
\end{cases}
\]
Similarly, we have $b_{0}\in J_0$, where $J_0$ is an interval of length $\mu_1 J_0=8c_{11}\cdot\max\{1,|b_2|\mu_1 I_1\}$ centered at $b_{0,0}$, and the number of possible values for $b_0$ does not exceed the measure of the interval $J_0$.

The following estimate
\begin{equation}\label{eq16}
\# (b_1,b_0)\leq \mu_1 J_1\cdot \mu_1 J_0=
\begin{cases}
2^7\varepsilon^{-1}c_{11}^2\cdot|b_2|^2\mu_2\Pi,\quad |b_2|\ge \left(\mu_1 I_1\right)^{-1},\\
2^7\varepsilon^{-1}c_{11}^2\cdot|b_2|\mu_1 I_2,\quad \left(\mu_1 I_2\right)^{-1}\leq|b_2|\leq \left(\mu_1 I_1\right)^{-1},\\
2^7\varepsilon^{-1}c_{11}^2,\quad |b_2|\leq \left(\mu_1 I_2\right)^{-1},
\end{cases}
\end{equation}
holds for a fixed value of the coefficient $b_2$.

Let us use the estimates (\ref{eq12}) and (\ref{eq16}) to consider the following three cases.

{\it Case 1}: $\left(\mu_1 I_1\right)^{-1}\leq |b_2|\leq  4\delta_2\varepsilon^{-1}Q$.

In this case, the first estimate of (\ref{eq16}) holds, and we have
\[
\mu_2 L_{2,2}\leq 2^{11}\varepsilon^{-3}c_{11}^2h_2^2\cdot Q^{-1}\mu_2\Pi\cdot 4\delta_2\varepsilon^{-1}Q<\textstyle\frac{1}{12r}\mu_2\Pi,
\]
for $\delta_1 < 2^{-17}r^{-1}\varepsilon^{4}c_{11}^{-2}h_2^{-2}$.

{\it Case 2}: $\left(\mu_1 I_2\right)^{-1}\leq |b_2|\leq  \left(\mu_1 I_1\right)^{-1}$.

Then the second estimate of (\ref{eq16}) holds, and we have
\[
\mu_2 L_{2,2}\ll Q^{-1}\mu_1 I_2\sum\limits_{\left(\mu_1 I_2\right)^{-1}\leq |b_2|\leq  \left(\mu_1 I_1\right)^{-1}}{|b_2|^{-1}}\ll Q^{-1}\ln Q\cdot\mu_1 I_2.
\]
Consequently, for $\varepsilon_1 = \textstyle\frac{1-s_1}{2}$ and $Q>Q_0$ we obtain
\[
\mu_2 L_{2,2}\ll Q^{-1+\varepsilon_1}\mu_1 I_2\ll Q^{-\varepsilon_1}\mu_2\Pi \leq\textstyle\frac{1}{12r}\mu_2\Pi.
\]

{\it Case 3}: $1\leq |b_2|\leq  \left(\mu_1 I_2\right)^{-1}$.

In this case, the third estimate of (\ref{eq16}) holds, leading to
\[
\mu_2 L_{2,2}\leq 2^{11}\varepsilon^{-3}c_{11}^2h_2^2\cdot Q^{-1}\sum\limits_{1\leq |b_2|\leq  \left(\mu_1 I_2\right)^{-1}}|b_2|^{-2}\leq\textstyle\frac{1}{12r}\mu_2\Pi,
\]
for $c_{1,1}c_{1,2} > 2^{12}r\pi^2c_{11}^2\varepsilon^{-3}h_2^2$.
\end{Proof}
 
\subsection{The induction step: reducing the degree of the polynomial.}\label{sec_2}

Let us return to the proof of Lemma \ref{lm5}. For $|P'(\alpha_1)|\in T_{1,j_1}$ and $|P'(\alpha_2)|\in T_{2,j_2}$, we have the following system of inequalities:
\begin{equation}\label{eq17}
\begin{cases}
|P(x_i)|< h_n\cdot Q^{-v_i},\quad v_i > 0,\\
d_{j_i,i}Q^{t_{j_i,i}}\leq |P'(\alpha_i)|<d_{j_i-1,i}Q^{t_{j_i-1,i}},\\
v_1+v_2=n-1,\quad i=1,2.
\end{cases}
\end{equation}

Without loss of generality, assume that $j_1\leq j_2$. Denote by $L_{j_1,j_2}$ the set of points $(x_1,x_2)\in\Pi$ such that the system of inequalities (\ref{eq17}) has a solution in polynomials $P\in\Pl_n(Q)$. By Lemma \ref{lm1}, it follows that $L_{j_1,j_2}$ is contained in a union $\bigcup\limits_{P\in\Pl_n(Q)}{\sigma_P}$, where
\begin{equation}\label{eq24}
\sigma_P=
\left\{(x_1,x_2)\in\Pi:\quad
|x_i-\alpha_i|<2^{n-1}h_n\cdot Q^{-v_i}|P'(\alpha_i)|^{-1},i=1,2\right\}.
\end{equation}
It means that the following estimate for $\mu_2 L_{j_1,j_2}$ holds:
\[
\mu_2 L_{j_1,j_2}\leq\mu_2\bigcup\limits_{P\in\Pl_n(Q)}{\sigma_P}\leq\sum\limits_{P\in\Pl_n(Q)}{\mu_2\sigma_P}.
\]

Together with the sets $\sigma_P$ consider the following expanded sets
\begin{equation}\label{eq22}
\sigma'_P= \sigma'_{P,1}\times\sigma'_{P,2}=
\left\{(x_1,x_2)\in\Pi:\quad |x_i-\alpha_i|<c_{12}Q^{-\gamma_{j_2,i}}|P'(\alpha_i)|^{-1},i=1,2\right\}.
\end{equation}
where $\gamma_{j_2,i}=\frac{(j_2-1)v_i}{n-1}$. Simple calculations show that the measure of the set $\sigma'_P$ is smaller than the measure of the rectangle $\Pi$ for $Q>Q_0$:
\[
\mu_2\sigma'_P\leq 4c_{12}^2\cdot Q^{1-j_2}Q^{-t_{1,j_1}-t_{2,j_2}}<4c_{12}^2\cdot Q^{-\frac{j_2+1}{2}}<\mu_2\Pi.
\]
Using (\ref{eq24}) and (\ref{eq22}), we find that the measures $\mu_2\sigma_P$ and $\mu_2\sigma'_P$ are connected as follows:
\begin{equation}\label{eq18}
\mu_2\sigma_P\leq 2^{2n-2}h_n^2c_{12}^{-2}\cdot Q^{-n+j_2}\mu_2\sigma'_P.
\end{equation}

Fix the vector $\mathbf{b}_{j_2}=\left(a_n,\ldots,a_{j_2+1}\right)$, where $a_n,\ldots,a_{j_2+1}$ are the coefficients of the polynomial $P\in \Pl_n(Q)$. Denote by $\Pl_n(\mathbf{b}_{j_2})\subset\Pl_n(Q)$ a subclass of polynomials with the same vector of coefficients $\mathbf{b}_{j_2}$. The number of subclasses $\Pl_n(\mathbf{b}_{j_2})$ is equal to the number of vectors $\mathbf{b}_{j_2}$ which can be estimated as follows:
\begin{equation}\label{eq19}
\#\left\{\mathbf{b}_{j_2}\right\}=(2Q+1)^{n-j_2}<2^{2n}Q^{n-j_2}.
\end{equation}

We are going to apply Sprind\v{z}uk's method of essential and non-essential sets \cite{Sprindzuk67}. A set $\sigma'_{P_1}$, $P_1\in\Pl_n(\mathbf{b}_{j_2})$ is called 
{\it essential} if for every  $\sigma'_{P_2}$,
$P_2\in\Pl_n(\mathbf{b}_{j_2})$, $P_2\neq P_1$, the inequality
\begin{equation}
\mu_2\left(\sigma'_{P_1}\cap\sigma'_{P_2}\right) < \textstyle\frac12 \mu_2\sigma'_{P_1},
\end{equation}
is satisfied.
Otherwise, $\sigma'_{P_1}$ is called {\it non-essential}.

{\it The case of essential sets.} For essential sets, we have the following estimate:
\begin{equation}\label{eq20}
\sum\limits_{\substack{P\in\Pl_n(\mathbf{b}_{j_2})\\ \sigma'_{P}-\text{essential}}}\mu_2\sigma'_{P}\leq 4\mu_2\Pi.
\end{equation}
Then from (\ref{eq18}), (\ref{eq19}) and (\ref{eq20}) we can write
\begin{equation}\label{eq21}
\sum\limits_{\mathbf{b}_{j_2}}\sum\limits_{\substack{P\in\Pl_n(\mathbf{b}_{j_2})\\ \sigma'_{P}-\text{essential}}}\mu_2\sigma_P\leq 2^{4n-2}h_n^2c_{12}^{-2}\sum\limits_{\substack{P\in\Pl_n(\mathbf{b}_{j_2})\\ \sigma'_{P}-\text{essential}}}\mu_2\sigma'_P< \textstyle\frac{1}{24r}\mu_2\Pi,
\end{equation}
for $c_{12}=2^{2n+3}r^{1/2}h_n$.

{\it The case of non-essential sets.} If a set $\sigma'_{P_1}$ is non-essential, then there exists a set $\sigma'_{P_2}$ such that $\mu_2\left(\sigma'_{P_1}\cap\sigma'_{P_2}\right) > \textstyle\frac12 \mu_2\sigma'_{P_1}$. Consider the polynomial $R= P_2-P_1$, $\deg R\leq j_2$, $H(R)\leq 2Q$, on the set $\left(\sigma'_{P_1}\cap\sigma'_{P_2}\right)$. Let us estimate the values $|R(x_i)|$ and $|R'(x_i)|$, $i,j=1,2$.

Let us write Taylor expansions of the polynomials $P_1$ and $P_2$ in the interval $\sigma'_{P_1,i}\cap\sigma'_{P_2,i}$, $i=1,2$:
\[
P_j(x_i)=P_j'(\alpha_{j,i})(x_i-\alpha_{j,i})+\ldots+\textstyle\frac{1}{n!}\cdot P_j^{(n)}(\alpha_{j,i})(x_i-\alpha_{j,i})^n,
\]
where $\alpha_{j,i}\in\sigma'_{P_j,i}$. From the estimate (\ref{eq22}), we have:
\[
|P_j'(\alpha_{j,i})(x_i-\alpha_{j,i})|\leq c_{12}Q^{-\gamma_{j_2,i}},
\]
\[
\left|\textstyle\frac{1}{k!}P_j^{(k)}(\alpha_{j,i})(x_i-\alpha_{j,i})^k\right|\leq c_{13,k} Q^{1-k\gamma_{j_2,i}-kt_{j_2,i}}\leq c_{13,k} Q^{1-\frac{k}{2}+\frac{k}{2}\gamma_{j_2,i}-k\gamma_{j_2,i}}\leq c_{13,k}Q^{-\gamma_{j_2,i}},
\]
for $k\ge 2$ and $Q>Q_0$.

Thus, the estimate $|R(x_i)|<|P_1(x_i)|+|P_2(x_i)|<c_{13}\cdot Q^{-\gamma_{j_2,i}}$ holds. From Lemma \ref{lm2} it follows that for every point $(x_1,x_2)\in\sigma'_{P_1}$, the inequalities
\[
|R(x_i)|<c_{14}\cdot Q^{-\gamma_{j_2,i}},\quad i=1,2,
\]
are satisfied.

Now let us write Taylor expansions of the polynomials $P'_1$ and $P'_2$ in the interval $\sigma'_{P_1,i}\cap\sigma'_{P_2,i}$, $j,i=1,2$:
\[
P'_j(x_i)=P'_j(\alpha_{j,i})+\ldots+\textstyle\frac{1}{(n-1)!} P_j^{(n)}(\alpha_{j,i})(x_i-\alpha_{j,i})^{n-1},
\]
where $\alpha_{j,i}\in\sigma'_{P_j,i}$. From the estimate (\ref{eq22}), we have:
\[
\left|\textstyle\frac{1}{(k-1)!}P_j^{(k)}(\alpha_i)(x_i-\alpha_i)^{k-1}\right|\leq c_{15,k}Q^{1+(k-1)\left(\frac{\gamma_{j_2,i}}{2}-\gamma_{j_2,i}-\frac12\right)}\leq c_{15,k}|P'(\alpha_i)|
\]
for $Q>Q_0$. Thus, we obtain $|R'(x_i)|\leq|P'_1(x_i)|+|P'_2(x_i)|\leq c_{15}|P'(\alpha_i)|$.  From Lemma \ref{lm2} it follows that for a sufficiently large $Q>Q_0$ the following inequalities hold:
\[
\min\limits_i\left\{|R'(x_i)|\right\}\leq c_{16}\min\limits_i\left\{|P'(\alpha_i)|\right\}\leq \begin{cases}
2c_{16}\delta_n Q,\quad j_1=j_2=2,\\
c_{16}Q^{\frac12},\quad j_1\neq 2 \text{ or }j_2\neq 2,
\end{cases}
\]
for every point $(x_1,x_2)\in\sigma'_{P_1}$. Thus, the measure of $L_{j_1,j_2}$ for non-essential sets does not exceed the respective measure for the system
\begin{equation}\label{eq23}
\begin{cases}
|R(x_i)|< h_{j_2}Q_1^{-\gamma_{j_2,i}},\quad \gamma_{j_2,i} > 0,\\
\min\limits_i\left\{|R'(x_i)|\right\}<\delta_{j_2}Q_1,\\
\gamma_{1,j_2}+\gamma_{2,j_2}=j_2-1,\quad i = 1,2,
\end{cases}
\end{equation}
where $Q_1=\min\limits_i\left\{(h_{j_2}/c_{14})^{1/\gamma_{j_2,i}}\right\}\cdot Q$ and $\delta_{j_2}=2c_{16}\cdot\left(\min\limits_i\left\{(h_{j_2}/c_{14})^{1/\gamma_{j_2,i}}\right\}\right)^{-1}\cdot\delta_n$. 

It should be mentioned that if polynomial $R(t)=a_1t-a_0$ is linear, then by Lemma \ref{lm1} we obtain:
\[
\left|x_i-\textstyle\frac{a_0}{a_1}\right|\ll Q_1^{-\gamma_{j_2,i}}<\textstyle\frac{\varepsilon}{4},\quad i=1,2
\]
for $Q_1>Q_0$. Hence, we immediately have $|x_1-x_2|<\varepsilon$ which contradicts to condition 2 for polynomial $\Pi$. Thus, $\deg R\ge 2$ and we can use induction. Since $j_2<n$, by the induction hypothesis the measure of solutions of the system (\ref{eq23}) is bounded from above by $\frac{1}{24r}\mu_2\Pi$ for $\delta_{j_2} \leq \delta_0$ and $Q_1>Q_0$. Thus, 
\[
\sum\limits_{\mathbf{b}_{j_2}}\sum\limits_{\substack{P\in\Pl_n(\mathbf{b}_{j_2})\\ \sigma'_{P}-\text{non-essential}}}\mu_2\sigma_P\leq\sum\limits_{\mathbf{b}_{j_2}}\sum\limits_{\substack{P\in\Pl_n(\mathbf{b}_{j_2})\\ \sigma'_{P}-\text{non-essential}}}\mu_2\sigma'_P\leq \textstyle\frac{1}{24r}\mu_2\Pi
\]
and together with the estimate (\ref{eq21}), this implies that
\[
\mu_2 L_{j_1,j_2}\leq \sum\limits_{\mathbf{b}_{j_2}}\sum\limits_{\substack{P\in\Pl_n(\mathbf{b}_{j_2})\\ \sigma'_{P}-\text{essential}}}\mu_2\sigma_P + \sum\limits_{\mathbf{b}_{j_2}}\sum\limits_{\substack{P\in\Pl_n(\mathbf{b}_{j_2})\\ \sigma'_{P}-\text{non-essential}}}\mu_2\sigma_P\leq \textstyle\frac{1}{12r}\mu_2\Pi.
\]

\subsection{The case of sub-intervals $T_{1,n}$ and $T_{2,n}$}\label{sec_3}

For $|P'(\alpha_1)|\in T_{1,n}$ and $|P'(\alpha_2)|\in T_{2,n}$ we have the following system of inequalities:
\begin{equation}\label{eq25}
\begin{cases}
|P(x_i)|< h_n\cdot Q^{-v_i},\quad v_i>0,\\
c_{10}Q^{\frac12-\frac{v_i}{2}}\leq |P'(\alpha_i)|<Q^{\frac12-\frac{v_i}{2}+\frac{v_i}{2(n-1)}},\\
v_1+v_2=n-1,\quad i=1,2.
\end{cases}
\end{equation}

By Lemma \ref{lm1}, the set $L_{n,n}$ of solutions of the system (\ref{eq25}) is contained in a union $\bigcup\limits_{P\in\Pl_n(Q)}{\sigma_P}$, where
\begin{equation}\label{eq26}
\sigma_P=\left\{(x_1,x_2)\in\Pi:\quad
|x_i-\alpha_i|\leq 2^{n-1}h_nc_{10}^{-1} Q^{-\frac{v_i+1}{2}},i=1,2\right\}.
\end{equation}

This leads to the following estimate for $\mu_2 L_{n,n}$:
\[
\mu_2 L_{n,n}\leq\mu_2\bigcup\limits_{P\in\Pl_n(Q)}{\sigma_P}\leq\sum\limits_{P\in\Pl_n(Q)}{\mu_2\sigma_P}.
\]

In this case we can not apply induction since the degree of the polynomial can not be reduced. Let us use a different method to estimate the measure $\mu_2 L_{n,n}$.

Cover the rectangle $\Pi$ by a system of disjoint rectangles $\Pi_k=J_{1,k}\times J_{2,k}$, where $\mu_1 J_{i,k}=Q^{-\frac{v_i+1}{2}+\varepsilon_{2,i}}$, $i=1,2$, such that $\Pi\subset\bigcup\limits_k{\Pi_k}$ and $\Pi_k\cap\Pi\neq\emptyset$. Thus, the number of rectangles $\Pi_k$ can be estimated as follows:
\begin{equation}\label{eq27}
2\max\left\{\frac{\mu_1 I_1}{\mu_1 J_{1,k}}, 1\right\}\cdot 2\max\left\{\frac{\mu_1 I_2}{\mu_1 J_{2,k}}, 1\right\}=
\begin{cases}
4Q^{\frac{n+1}{2}-\varepsilon_{2,1}-\varepsilon_{2,2}}\mu_2\Pi,\quad s_i<\frac{v_i+1}{2},\\
4Q^{\frac{v_1+1}{2}-\varepsilon_{2,1}}\mu_1 I_1,\quad s_1<\frac{v_1+1}{2}, s_2\ge\frac{v_2+1}{2},\\ 
4Q^{\frac{v_2+1}{2}-\varepsilon_{2,2}}\mu_1 I_2,\quad s_1\ge\frac{v_1+1}{2}, s_2<\frac{v_2+1}{2}. 
\end{cases}
\end{equation}

We are going to say that a polynomial $P$ belongs to $\Pi_k$ if there is a point $(x_1,x_2)\in\Pi_k$ such that the inequalities (\ref{eq25}) are satisfied. 

Now let us prove that there is no rectangle $\Pi_k$ containing two or more irreducible polynomials $P\in\Pl_n(Q)$. Assume the converse: let $P_1, P_2\in\Pi_k$ be irreducible polynomials and let the inequalities (\ref{eq25}) hold for each polynomial $P_j$ at a point $(x_{j,1},x_{j,2})\in\Pi_k$, $j=1,2$. Thus, for $Q>Q_0$ and for every point $(x_1,x_2)\in\Pi_k$, the estimates
\begin{equation}\label{eq28}
|x_i-\alpha_{j,i}|\leq|x_i-x_{j,i}|+|x_{j,i}-\alpha_{j,i}|\leq 2Q^{-\frac{v_i+1}{2}+\varepsilon_{2,i}},
\end{equation}
are satisfied, where $x_{j,i}\in S(\alpha_{j,i})$.

Let us estimate the values $|P_j(x_i)|$, $i,j=1,2$ where $(x_1,x_2)\in\Pi_k$. Let us write Taylor expansions of $P_j$ in the interval $J_{i,k}$:
\[
P_j(x_i)=P'_j(\alpha_{j,i})(x_i-\alpha_{j,i})+\ldots+\textstyle\frac{1}{n!}\cdot P^{(n)}_j(\alpha_{j,i})(x_i-\alpha_{j,i})^{n}.
\]
From estimates (\ref{eq25}) and (\ref{eq28}) we obtain that
\[
\left|P'_j(\alpha_{j,i})(x_i-\alpha_{j,i})\right|\ll Q^{-v_i+\frac{v_i}{2(n-1)}+\varepsilon_{2,i}},
\]
\[
\left|\textstyle\frac{1}{k!}\cdot P^{(k)}_j(\alpha_{j,i})(x_i-\alpha_{j,i})^k\right|\ll Q^{1-\frac{k}{2}-\frac{kv_i}{2}+k\varepsilon_{2,i}}\ll Q^{-v_i+\frac{v_i}{2(n-1)}+\varepsilon_{2,i}}
\]
for $\varepsilon_{2,i}<\textstyle\frac{v_i}{2(n-1)(k-1)}$ and $Q>Q_0$.

Then we can write the following estimate:
\begin{equation}\label{eq0050}
|P_j(x_i)|\ll Q^{-v_i+\frac{v_i}{2(n-1)}+\varepsilon_{2,i}}<Q^{-v_i+\frac{v_i}{2(n-1)}+\varepsilon_{2,i}+\varepsilon_3},
\end{equation}
where $\varepsilon_{2,i}<\textstyle\frac{v_i}{2(n-1)^2}$.

From Lemma \ref{lm3} for $\eta_i = \textstyle\frac{v_i+1}{2}-\varepsilon_{2,i}$ and $\tau_i=v_i-\textstyle\frac{v_i}{2(n-1)}-\varepsilon_{2,i}-\varepsilon_3$, $i=1,2$, we have
\[
\tau_1+\tau_2+2=(n-1)-\textstyle\frac12-\varepsilon_{2,1}-\varepsilon_{2,2}+2-2\varepsilon_3=n+\textstyle\frac12-\varepsilon_{2,1}-\varepsilon_{2,2}-2\varepsilon_3,
\]
\[
2(\tau_i+1-\eta_i)= 2\left(v_i-\textstyle\frac{v_i}{2(n-1)}-\varepsilon_{2,i}-\varepsilon_3+1-\textstyle\frac{v_i+1}{2}+\varepsilon_{2,i}\right)=v_i+1-\textstyle\frac{v_i}{n-1}-2\varepsilon_3.
\]
Substitution of this expressions into (\ref{eq0}) leads to the inequality
\[
\tau_1+\tau_2+2+2(\tau_1+1-\eta_1)+2(\tau_2+1-\eta_2)=2n+\textstyle\frac12-\varepsilon_{2,1}-\varepsilon_{2,2}-6\varepsilon_3\ge2n+\textstyle\frac{1}{8}
\]
for $\varepsilon_{2,i}=\frac{v_i}{4(n-1)^2}$, $\varepsilon_3=\textstyle\frac{1}{48}$. This contradict to Lemma \ref{lm3} with $\delta=\frac18$.

Hence, every rectangle $\Pi_k$ contains at most one polynomial $P\in\Pl_n(Q)$. In this case, we have the following estimate for the measure of the set $L_{n,n}$:
\[
\mu_2 L_{n,n}\leq\sum\limits_{\Pi_k}{\mu_2\sigma_P},
\]
and together with the estimates (\ref{eq26}) and (\ref{eq27}) this leads to
\[
\mu_2 L_{n,n}\ll Q^{-\varepsilon_{2,1}-\varepsilon_{2,2}}\mu_2\Pi<\textstyle\frac{1}{12r}\mu_2\Pi
\]
for $Q>Q_0$ and $s_i<\frac{v_i+1}{2}$, $i=1,2$. If $s_i \ge \frac{v_i+1}{2}$, then we obtain the estimate 
\[
\mu_2 L_{n,n}\leq\sum\limits_{P\in\Pl_n(Q)}{\mu_2\sigma_P}\ll Q^{-\varepsilon_{2,i}}\mu_1 I_1\mu_1 I_2<\textstyle\frac{1}{12r}\mu_2\Pi
\]
for $Q>Q_0$.

\subsection{The case of a small derivative}\label{sec_4}

Let us discuss a situation where $|P'(x_i)|\leq 2c_{10}Q^{\frac12-\frac{v_i}{2}}$, $i=1,2$. In this case, we can show that $|P'(\alpha_i)|\leq 2^{n-1}c_{10}Q^{\frac12-\frac{v_i}{2}}$, where $x_i\in S(\alpha_i)$.

Indeed, let $|P'(\alpha_i)|> 2^{n-1}c_{10}Q^{\frac12-\frac{v_i}{2}}$. Let us write a Taylor expansions of the polynomial $P'$:
\[
P'(x_i)=P'(\alpha_i)+P''(\alpha_i)(x_i-\alpha_i)+\ldots+\textstyle\frac{1}{(n-1)!} P^{(n)}(\alpha_i)(x_i-\alpha_i)^{n-1}.
\]
Using our assumption and repeating analogous computations to those from the beginning of the proof of Lemma \ref{lm5} (see page 5) we have:
\[
\left|P''(\alpha_i)(x_i-\alpha_i)+\ldots+\textstyle\frac{1}{(n-1)!} P^{(n)}(\alpha_i)(x_i-\alpha_i)^{n-1}\right|\leq c_{10}Q^{\frac12-\frac{v_i}{2}}.
\]
This leads to the following upper bound for $|P'(\alpha_i)|$:
\[
|P'(\alpha_i)|\leq 3c_{10}Q^{\frac12-\frac{v_i}{2}},
\]
which contradicts our assumption for $n\ge 3$.

Now let $L_{n+1,n+1}\subset\Pi$ be the set of points satisfying the system
\begin{equation}\label{eq37}
\begin{cases}
|P(x_i)|< h_n Q^{-v_i},\quad v_i >0,\\
|P'(\alpha_i)|<2^{n-1}c_{10}Q^{\frac12-\frac{v_i}{2}},\\
v_1+v_2=n-1,\quad i=1,2.
\end{cases}
\end{equation}

The polynomials $P\in \mathcal{P}_n(Q)$ satisfying (\ref{eq37}) are going to be classified according to the distribution of their roots and the size of the leading coefficient $|a_m|$. This classification was introduced by Sprind\v{z}uk \cite{Sprindzuk67}.

For every polynomial $P\in \Pl_n(Q)$ of degree $3\leq m\leq n$ we define numbers $\rho_{1,j}$ and $\rho_{2,j}$, $2\le j\le m$, as solutions of equations
\[
|\alpha_{1,1}-\alpha_{1,j}| = Q^{-\rho_{1,j}},\quad|\alpha_{2,1}-\alpha_{2,j}| = Q^{-\rho_{2,j}}.
\]
Let us also define the vectors $\mathbf{k}_1=(k_{1,2},\ldots,k_{1,m})$ and $\mathbf{k}_2=(k_{2,2},\ldots,k_{2,m})$ with integer coefficients as solutions of the inequalities
\[
k_{i,j}\varepsilon_4-\varepsilon_4 \le \rho_{i,j} < k_{i,j}\varepsilon_4,\quad i=1,2, j=\overline{2,m},
\]
where $\varepsilon_4 > 0$ is some small constant.

Denote by $\Pl_m(Q, \mathbf{k}_1, \mathbf{k}_2, u)\subset\Pl_n(Q)$ a subclass of polynomials with the same pair of vectors $(\mathbf{k}_1, \mathbf{k}_2)$ and the following bounds on leading coefficients: $Q^{u}\leq|a_m|<Q^{u+\varepsilon_4}$, where $u\in \Z\cdot \varepsilon_4$. Since $1\leq |a_m|\leq Q$, the following estimate holds for $u$: $0\leq u\leq 1-\varepsilon_4$. The roots of the polynomial $P$ are bounded, and we can write $Q\gg |\alpha_{j_1}-\alpha_{j_2}| \gg H^{-m+1} \gg Q^{-m+1}$, which leads to the estimates $-\textstyle\frac{1}{\varepsilon_4}\le k_{i,j} \le \textstyle\frac{m-1}{\varepsilon_4}+1$. Thus, an integer vector $\mathbf{k}_i = (k_{i,2},\ldots,k_{i,m})$ can take at most $\left(\textstyle\frac{m}{\varepsilon_4}+1\right)^{m-1}$ values, the number of subclasses $\Pl_m(Q, \mathbf{k}_1, \mathbf{k}_2, u)$ can be estimated as follows:
\begin{equation}\label{eq29}
\# \{m, \mathbf{k}_1, \mathbf{k}_2, l\}\leq nc_{16}^2c_{17},
\end{equation}
where $c_{16}=\sum\limits_{i=2}^n{\left(\textstyle\frac{i}{\varepsilon_4}+1\right)^{i-1}}$, $c_{17}=\varepsilon_4^{-1}+1$.

Let $p_{i,j}$, $i=1,2$, $j=\overline{1,m}$ be defined as follows:
\begin{equation}\label{eq34}
\begin{cases}
p_{i,j} = (k_{i,j+1} + \ldots + k_{i,m})\cdot\varepsilon_4, \qquad 1\le j \le m-1,\\
p_{i,j}= 0,\quad j=m.
\end{cases}
\end{equation}

For a polynomial $P\in \Pl_m(Q, \mathbf{k}_1, \mathbf{k}_2, u)$, we can write the following estimates for its derivatives at the root $\alpha_i$:
\begin{equation}\label{eq30}
\begin{aligned}
Q^{u-p_{i,1}}\leq |P'(\alpha_i)|&=|a_m||\alpha_{i,1}-\alpha_{i,2}|\ldots|\alpha_{i,1}-\alpha_{i,m}|\leq Q^{u-p_{i,1}+(m+1)\varepsilon_4},\\
|P^{(j)}(\alpha_i)|&\ll |a_m||\alpha_{i,1}-\alpha_{i,j+1}|\ldots|\alpha_{i,1}-\alpha_{i,m}|\ll Q^{u-p_{i,j}+(m+1)\varepsilon_4},\quad j=\overline{2,m}.
\end{aligned}
\end{equation}

Consider polynomials which solve the system (\ref{eq37}). We can assume that the following inequalities hold:
\[
\begin{aligned}
Q^{u-p_{1,1}}&\leq |P'(\alpha_1)|\ll Q^{\frac12-\frac{v_1}{2}},\\
Q^{u-p_{2,1}}&\leq |P'(\alpha_2)|\ll Q^{\frac12-\frac{v_2}{2}},
\end{aligned}
\]
which leads to the inequalities
\begin{equation}\label{eq38}
p_{1,1}>u+\frac{v_1-1}{2},\quad p_{2,1}>u+\frac{v_2-1}{2}.
\end{equation}

Now let us obtain an estimate for the measure of the set $L_{n+1,n+1}$. From Lemma \ref{lm1} it follows that this set is contained in a union $\bigcup\limits_{m,\mathbf{k}_1,\mathbf{k}_2,u}{\bigcup\limits_{P\in\Pl_m(Q,\mathbf{k}_1, \mathbf{k}_2, u)}{\sigma_P}}$, where
\[
\sigma_P=\left\{(x_1,x_2)\in\Pi:\quad
|x_i-\alpha_i| \le \min\limits_{1\le j\le m} \left(2^{m-j} \textstyle\frac{h_n\cdot Q^{-v_i}}{|P'(\alpha_{i,1})|}\cdot|\alpha_{i,1}-\alpha_{i,2}|\ldots|\alpha_{i,1}-\alpha_{i,j}| \right)^{1/j}\right\}.
\]

This, together with the previous notation (\ref{eq34}) and the estimates (\ref{eq30}), yields the formula
\begin{equation}\label{eq035}
\sigma_P=\left\{(x_1,x_2)\in\Pi:\quad
|x_i-\alpha_i| \leq \textstyle\frac12\cdot \min\limits_{1\le j\le m} \left((2^mh_n)^{1/j}\cdot Q^{\frac{-u-v_i+p_{i,j}}{j}}\right),i=1,2\right\}
\end{equation}
for $P\in \Pl_m(Q, \mathbf{k}_1, \mathbf{k}_2, u)$. If the inequalities
\begin{equation}\label{eq33}
(2^mh_n)^{1/m_i}\cdot Q^{\frac{-u-v_i+p_{i,m_i}}{m_i}}\leq (2^mh_n)^{1/k}\cdot Q^{\frac{-u-v_i+p_{i,k}}{k}},\quad 1\leq k\leq m, i=1,2,
\end{equation}
are satisfied, then the numbers $j=m_1$ and $j=m_2$ provide the best estimates for the roots $\alpha_1$ and $\alpha_2$ respectively, and the inequalities
\begin{equation}\label{eq32}
\sigma_P=\left\{(x_1,x_2)\in\Pi:\quad
|x_i-\alpha_i| \leq \textstyle\frac12\cdot(2^mh_n)^{1/m_i}\cdot Q^{\frac{-u-v_i+p_{i,m_i}}{m_i}},i=1,2\right\}
\end{equation}
hold.

Let us cover the rectangle $\Pi$ by a system of disjoint rectangles $\Pi_{m_1,m_2}=J_{m_1}\times J_{m_2}$, where $\mu_1 J_{m_i}=Q^{-\frac{u+v_i-p_{i,m_i}}{m_i}+\varepsilon_5}$, such that $\Pi\subset\bigcup\limits_k{\Pi_{m_1,m_2}}$ and $\Pi_{m_1,m_2}\cap\Pi\neq\emptyset$. The number of rectangles $\Pi_{m_1,m_2}$ can be estimated as follows:
\begin{equation}\label{eq35}
\# \Pi_{m_1,m_2}\leq 4\cdot Q^{\frac{u+v_1-p_{1,m_1}}{m_1}+\frac{u+v_2-p_{2,m_2}}{m_2}-2\varepsilon_5}\mu_2\Pi.
\end{equation}

Now let us show that there is no rectangle $\Pi_{m_1,m_2}$ containing two or more irreducible polynomials. Let $P_1,P_2\in\Pi_{m_1,m_2}$ be irreducible polynomials, and let the inequalities (\ref{eq37}) hold for polynomials $P_j$ at points $(x_{j,1},x_{j,2})\in\Pi_{m_1,m_2}$, $j=1,2$. Thus, estimates
\begin{equation}\label{eq36}
|x_i-\alpha_{j,i}|\leq|x_i-x_{j,i}|+|x_{j,i}-\alpha_{j,i}|\leq 2\cdot Q^{-\frac{u+v_i-p_{i,m_i}}{m_i}+\varepsilon_5}
\end{equation}
are satisfied for every point $(x_1,x_2)\in\Pi_{m_1,m_2}$ and for $Q>Q_0$, where $x_{j,i}\in S(\alpha_{j,i})$.

Let us estimate $|P_j(x_i)|$, where $(x_1,x_2)\in\Pi_{m_1,m_2}$. Let us write Taylor expansions of the polynomials $P_j$ in the interval $J_{m_i}$:
\[
P_j(x_i)=P'_j(\alpha_{j,i})(x_i-\alpha_{j,i})+\ldots +\textstyle\frac{1}{m!}\cdot P^{(m)}_j(\alpha_{j,i})(x_i-\alpha_{j,i})^m.
\]
By estimates (\ref{eq30}), (\ref{eq33}) and (\ref{eq36}) we have
\[
\left|\textstyle\frac{1}{k!}\cdot P^{(k)}_j(\alpha_{j,i})(x_i-\alpha_{j,i})^k\right|\ll Q^{-v_i+(m+1)\varepsilon_4+k\varepsilon_5}.
\]

This leads to the following estimates for $|P_j(x_i)|$: 
\begin{equation}\label{eq00050}
|P_j(x_i)|\ll Q^{-v_i+(m+1)\varepsilon_4+m\varepsilon_5}<Q^{-v_i+(m+1)(\varepsilon_4+\varepsilon_5)}.
\end{equation}

From Lemma \ref{lm3} with $\eta_i = \frac{u+v_i-p_{i,m_i}}{m_i}-\varepsilon_5$ and $\tau_i=v_i-(m+1)(\varepsilon_4+\varepsilon_5)$, where $i=1,2$ and $\varepsilon_4=\textstyle\frac{1}{12(m+1)}$, $\varepsilon_5=\textstyle\frac{1}{4(3m+1)}$, we obtain
\[
\tau_1+\tau_2+2=n+1-\textstyle\frac16-2(m+1)\varepsilon_5,
\]
\[
2(\tau_i+1-\eta_i)= 2v_i+2-2\cdot\textstyle\frac{u+v_i-p_{i,m_i}}{m_i}-\textstyle\frac16-2m\varepsilon_5.
\]
Let us estimate the expression $2(\tau_i+1-\eta_i)$ by applying the inequalities (\ref{eq38}):
\[
2(\tau_i+1-\eta_i)\ge
\begin{cases}
v_i+2-u+\textstyle\frac{2p_{i,m_i}}{m}-\textstyle\frac16-2m\varepsilon_5,\quad m_i\ge 2,\\
v_i+1-\textstyle\frac16-2m\varepsilon_5,\quad m_i=1,
\end{cases}
\ge v_i+1-\textstyle\frac16-2m\varepsilon_5.
\]
Substituting this expressions into (\ref{eq0}) yields
\[
\tau_1+\tau_2+2+2(\tau_1+1-\eta_1)+2(\tau_2+1-\eta_2)=2n+\textstyle\frac32-(6m+2)\varepsilon_5>2n+\textstyle\frac12,
\]
which contradicts to Lemma \ref{lm3} with $\delta=\frac12$.

This means that every rectangle $\Pi_{m_1,m_2}$ contains at most one polynomial $P\in\Pl_m(Q, \mathbf{k}_1, \mathbf{k}_2, u)$. Thus, the measure of solutions of the system (\ref{eq37}) can be estimated as follows:
\[
\mu_2 L_{n+1,n+1}\leq\sum\limits_{m,\mathbf{k}_1, \mathbf{k}_2, u}{\sum\limits_{P\in\Pl_m(Q, \mathbf{k}_1, \mathbf{k}_2, u)}{\mu_2\sigma_P}}\leq\sum\limits_{m,\mathbf{k}_1, \mathbf{k}_2, u}{\sum\limits_{\Pi_{m_1,m_2}}{\mu_2\sigma_P}}.
\]
Thus, by estimates (\ref{eq29}), (\ref{eq32}) and (\ref{eq35}), we can obtain the inequality
\[
\mu_2 L_{n+1,n+1}\ll Q^{-2\varepsilon_5}\cdot\mu_2\Pi<\textstyle\frac{1}{12r}\mu_2\Pi
\]
for $Q>Q_0$.

\subsection{Mixed cases}\label{sec_5}

{\bf The case of sub-intervals $T_{1,n}$, $T_{2,j}$ ($T_{1,j}$, $T_{2,n}$), $j=\overline{2, n-1}$}

Consider the system of inequalities
\begin{equation}\label{eq39}
\begin{cases}
|P(x_i)|< h_n\cdot Q^{-v_i},\quad v_i>0,\\
c_{10}Q^{\frac12-\frac{v_1}{2}}\leq|P'(\alpha_1)|<Q^{\frac12-\frac{v_1}{2}+\frac{v_1}{2(n-1)}},\\
Q^{\frac12-\frac{(j-1)v_2}{2(n-1)}}\leq|P'(\alpha_2)|<Q^{\frac12-\frac{v_2(j-2)}{2(n-1)}},\\
v_1+v_2=n-1,\quad i = 1,2.
\end{cases}
\end{equation}
Let $L_{n,j}$ be the set of solutions of the system (\ref{eq39}). In this case we need to consider two different sets. Let $L_{n,j}^1$ and $L_{n,j}^2$ be the sets of points $(x_1,x_2)\in\Pi$ such that there exists a polynomial $P\in\Pl_n(Q)$ satisfying the system (\ref{eq39}) under condition $c_{10}Q^{\frac12-\frac{v_1}{2}}\leq|P'(\alpha_1)|<Q^{\frac12-\frac{v_1}{2}+\frac{v_1}{4(n-1)}}$ and $Q^{\frac12-\frac{v_1}{2}+\frac{v_1}{4(n-1)}}\leq|P'(\alpha_1)|<Q^{\frac12-\frac{v_1}{2}+\frac{v_1}{2(n-1)}}$ respectively.

As in the case of small derivatives, we classify polynomials $P\in\Pl_n(Q)$ according to the distribution of their roots and the size of their leading coefficients. We will consider the subclasses of polynomials $\Pl_m(Q,\mathbf{k}_2,u)$ with the same vector $\mathbf{k}_2$ and the following bounds on leading coefficient: $Q^{u}<|a_m|<Q^{u+\varepsilon_4}$, where $0\leq u\leq 1-\varepsilon_4$, $0<\varepsilon_4<1$ and $u\in \Z\cdot \varepsilon_4$. Then
\begin{equation}\label{eq39_0}
\#\{m, \mathbf{k}_2,u\}\leq nc_{17}\cdot c_{16}.
\end{equation}

From Lemma \ref{lm1}, the set $L_{n,j}^g$, $g=1,2$ is contained in a union $\bigcup\limits_{m,\mathbf{k}_2,u}\bigcup\limits_{P\in\Pl_m(Q,\mathbf{k}_2,u)}{\sigma_P}$, where
\begin{equation}\label{eq40}
\sigma_P=\left\{(x_1,x_2)\in\Pi:\quad
\begin{array}{ll}
|x_1-\alpha_1|\leq 2^{m-1}h_n\max\{c_{10}^{-1},1\}\cdot Q^{-\frac{v_1}{2}-\frac12-\frac{v_1(g-1)}{4(n-1)}},\\
|x_2-\alpha_2|\leq 2^{m-1}h_n Q^{-v_2+p_{2,1}-u}
\end{array}
\right\}.
\end{equation}

Define the value $l=v_2-p_{2,1}+u-k_{2,2}\varepsilon_4$ and let us write $l$ as $l=[l]+\{l\}$, where $[l]$ is the integer part of $l$ and $\{l\}$ is the fractional part. Now let us cover the rectangle $\Pi$ by a system of disjoint rectangles $\Pi_k=J_{1,k}\times J_{2,k}$, where $\mu_1 J_{1,k}=Q^{-\frac{v_1}{2}-\frac12-\frac{v_1(g-1)}{4(n-1)}+\varepsilon_6}$ and $\mu_1 J_{2,k}= Q^{-k_{2,2}\varepsilon_4-\{l\}}$, such that $\Pi\subset\bigcup\limits_k{\Pi_k}$ and $\Pi_k\cap\Pi\neq\emptyset$.
The number of rectangles $\Pi_k\in\Pi$ can be estimated as
\begin{equation}\label{eq42}
\#\{\Pi_k\}\leq 4Q^{\frac{v_1}{2}+\frac12+\frac{v_1(g-1)}{4(n-1)}+k_{2,2}\varepsilon_4-\varepsilon_6+\{l\}}\mu_2\Pi.
\end{equation}

Assume that every rectangle $\Pi_k$ contains no more than $2^mQ^{[l]+\frac{\varepsilon_6}{2}}$ points $(\alpha_1,\alpha_2)$, where $\alpha_1,\alpha_2$ are the roots of polynomial $P\in\Pl_m(Q,\mathbf{k}_2,u)$. Then by inequalities (\ref{eq39_0}), (\ref{eq40}) and (\ref{eq42}) it follows that the measure of the set $L_{n,j}^g$ can be estimated as:
\begin{equation}\label{eq43}
\mu_2 L_{n,j}^g \leq 2^{3m+4}nc_{10}c_{16}c_{17}\cdot Q^{-v_2+p_{2,1}-u+k_{2,2}\varepsilon_4-\frac{\varepsilon_6}{2}+[l]+\{l\}}\mu_2\Pi\leq Q^{-\frac{\varepsilon_6}{2}}\mu_2\Pi\leq \textstyle\frac{1}{24r}\mu_2\Pi,
\end{equation}
where $Q>Q_0$.

Now assume that there exists a rectangle $\Pi_k$ containing more than $2^mQ^{[l]+\frac{\varepsilon_6}{2}}$ polynomials $P_j\in\Pl_m(Q,\mathbf{k}_2,u)$. From the Taylor expansions of polynomials $P_j$ in the interval $J_{2,k}$, the estimates (\ref{eq30}) and condition $(\alpha_{j,1},\alpha_{j,2})\in \Pi_k$ it follows that
\[
\left|\textstyle\frac{1}{k!}\cdot P_j^{(k)}(\alpha_{j,2})(x_2-\alpha_{j,2})^k\right|\ll Q^{u-p_{2,k}+(m+1)\varepsilon_4-k\cdot k_{2,2}\varepsilon_4-k\{l\}}< Q^{u-p_{2,1}-k_{2,2}\varepsilon_4-\{l\}+(m+1)\varepsilon_4},
\]
which allows us to write
\begin{equation}\label{eq43}
|P_j(x_2)|< Q^{u-p_{2,1}-k_{2,2}\varepsilon_4-\{l\}+(m+2)\varepsilon_4}.
\end{equation}
Similarly, repeating the calculations by analogy with Section \ref{sec_3} (see inequality (\ref{eq0050})), we have 
\begin{equation}\label{eq44}
|P_j(x_1)|< Q^{-v_1+\frac{v_1}{4(n-1)}+2\varepsilon_6}
\end{equation}
for $\varepsilon_6 <\textstyle\frac{v_1}{(n-1)^2}$.

By Dirichlet's principle we can find at least $\left[Q^{\frac{\varepsilon_6}{2}}\right]+1$ polynomials from $\Pl_m(Q,\mathbf{k}_2, u)$ contained in $\Pi_k$ such that their coefficients $a_m,\ldots, a_{m+1-[l]}$ coincide. Let us call them $P_1, \ldots, P_{\left[Q^{\frac{\varepsilon_6}{2}}\right]+1}$. If $[l]=0$, then we can simply ignore this step. Let us consider the differences $R_{i,j}=P_i-P_j$, $1\leq i< j\leq \left[Q^{\frac{\varepsilon_6}{2}}\right]+1$.

From the inequalities (\ref{eq43}) and (\ref{eq44}), we obtain that at every point of the rectangle $\Pi_k$ the polynomials $R_{i,j}$ satisfy
\begin{equation}\label{eq41}
\begin{cases}
|R_{i,j}(x_1)|< 2Q^{-v_1+\frac{v_1}{4(n-1)}+2\varepsilon_6},\quad |R_{i,j}(x_2)|< 2Q^{u-p_{2,1}-k_{2,2}\varepsilon_4-\{l\}+(m+2)\varepsilon_4},\\
\deg R_{i,j} \leq m-[l]=m-v_2+p_{2,1}-u+k_{2,2}\varepsilon_4+\{l\}.
\end{cases}
\end{equation}

Assume that among polynomials $R_{i,j}$ we can find at least two polynomials without common roots. Then we can apply Lemma \ref{lm3} with $\tau_1= v_1-\textstyle\frac{v_1}{4(n-1)}-2\varepsilon_{6}$, $\tau_2 = -u+p_{2,1}+k_{2,2}\varepsilon_4+\{l\}-(m+2)\varepsilon_4$, $\eta_1=\frac{v_1}{2}+\frac12+\frac{v_1(g-1)}{4(n-1)}-\varepsilon_6$, $\eta_2=k_{2,2}\varepsilon_4+\{l\}$, so that we have
\[
\tau_1+1=v_1+1-\textstyle\frac{v_1}{4(n-1)}-2\varepsilon_{6},\quad \tau_2+1=1-u+p_{2,1}+k_{2,2}\varepsilon_4+\{l\}-(m+2)\varepsilon_4,
\]
\[
2(\tau_1+1-\eta_1)= v_1+1-\textstyle\frac{gv_1}{2(n-1)}-2\varepsilon_6,\quad 2(\tau_2+1-\eta_2)=2-2u+2p_{2,1}-2(m+2)\varepsilon_4.
\]
Substituting these expressions into (\ref{eq0}) for $\varepsilon_4=\textstyle\frac{1-\{l\}}{9(m+2)}$ and $\varepsilon_6=\textstyle\frac{1-\{l\}}{12}$ yields
\begin{multline*}
\tau_1+\tau_2+2+2(\tau_1+1-\eta_1)+2(\tau_2+1-\eta_2)= 2v_1+5-\textstyle\frac{(2+g)v_1}{4(n-1)}+3p_{2,1}+k_{2,2}\varepsilon_4-3u+\\+\{l\}
-3(m+2)\varepsilon_4-4\varepsilon_{6}\ge 2n-2v_2+2p_{2,1}+2k_{2,2}\varepsilon_4-2u + (p_{2,1}-k_{2,2}\varepsilon_4)+\\+(1-u)+\{l\}
+2-\textstyle\frac{2(1-\{l\})}{3}-\textstyle\frac{v_1}{n-1}\ge 2(m-[l])-\textstyle\frac{\{l\}}{3} +\textstyle\frac13
\end{multline*}
This inequality contradict to Lemma \ref{lm3} for $\delta =\frac{1-\{l\}}{3}$.

The case when among polynomials $R_{i,j}$, $1\leq i< j\leq \left[Q^{\frac{\varepsilon_6}{2}}\right]+1$ we can not find two polynomials without common roots is considered in \cite{BernikGoetzeKukso14}.

Hence, we obtain
\[
\mu_2 L_{n,j}\leq \mu_2 L_{n,j}^1+\mu_2 L_{n,j}^2 \leq \textstyle\frac{1}{12r}\mu_2\Pi.
\]

{\bf The case where one derivative is small and the other derivative lies in the sub-interval $T_{2,j}$, $j=\overline{2, n}$ ($T_{2,j}$)}

Given the estimate for $|P'(\alpha_1)|$ obtained in Section \ref{sec_4} for $|P'(x_1)|\leq 2c_{10}Q^{\frac12-\frac{v_1}{2}}$ consider the system of inequalities
\begin{equation}\label{eq90}
\begin{cases}
|P(x_i)|< h_n\cdot Q^{-v_i},\quad v_i>0,\\
|P'(\alpha_1)|<2^{n-1}c_{10}Q^{\frac12-\frac{v_1}{2}},\\
Q^{\frac12-\frac{(j-1)v_2}{2(n-1)}}\leq |P'(\alpha_2)|<Q^{\frac12-\frac{(j-2)v_2}{2(n-1)}},\\
v_1+v_2=n-1,\quad i = 1,2.
\end{cases}
\end{equation}
Denote by $L_{n+1,j}$ the set of points $(x_1,x_2)\in\Pi$ such that there exists a polynomial $P\in\Pl_n(Q)$ satisfying the system (\ref{eq90}). Once again let us classify polynomials $P\in\Pl_n(Q)$ according to the distribution of their roots and the size of leading coefficients. We will consider the subclasses of polynomials $\Pl_m(Q,\mathbf{k}_1, \mathbf{k}_2,u)$ defined above.

From Lemma \ref{lm1} by analogy with Section \ref{sec_4} (see inequality (\ref{eq035})) we conclude that the set $L_{n+1,j}$ is contained in a union $\bigcup\limits_{m,\mathbf{k}_1,\mathbf{k}_2,u}\bigcup\limits_{P\in\Pl_m(Q,\mathbf{k}_1,\mathbf{k}_2,u)}{\sigma_P}$, where
\[
\sigma_P=\left\{(x_1,x_2)\in\Pi:\quad
\begin{array}{ll}
|x_1-\alpha_1|\leq \textstyle\frac12\min\limits_{1\le j\le m} \left((2^mh_n)^{1/j}\cdot Q^{\frac{-u-v_1+p_{1,j}}{j}}\right),\\
|x_2-\alpha_2|\leq 2^{m-1}h_n Q^{-u-v_2+p_{2,1}}
\end{array}
\right\}
\]
for $P\in \Pl_m(Q, \mathbf{k}_1, \mathbf{k}_2,u)$.

If the inequalities (\ref{eq33}) hold for $i=1$, then the estimate numbered as $j=m_1$ is optimal for the root $\alpha_1$, and we have
\begin{equation}\label{eq91}
\sigma_P=\left\{(x_1,x_2)\in\Pi:\quad
\begin{array}{ll}
|x_1-\alpha_1| \leq \textstyle\frac12\cdot(2^mh_n)^{1/m_1}\cdot Q^{\frac{-u-v_1+p_{1,m_1}}{m_1}},\\
|x_2-\alpha_2|\leq 2^{m-1}h_n Q^{-u-v_2+p_{2,1}}
\end{array}
\right\}.
\end{equation}

Define the value $l=v_2-p_{2,1}+u-k_{2,2}\varepsilon_4$ as in the previous case and let us cover the rectangle $\Pi$ by a system of disjoint rectangles $\Pi_k=J_{1,k}\times J_{2,k}$, where $\mu_1 J_{1,k}=Q^{-\frac{u+v_1-p_{1,m_1}}{m_1}+\varepsilon_7}$ and $\mu_1 J_{2,k}= Q^{-k_{2,2}\varepsilon_4-\{l\}}$, such that $\Pi\subset\bigcup\limits_k{\Pi_k}$ and $\Pi_k\cap\Pi\neq\emptyset$.
The number of rectangles $\Pi_k\in\Pi$ can be estimated as
\begin{equation}\label{eq92}
\#\{\Pi_k\}\leq 4Q^{\frac{u+v_1-p_{1,m_1}}{m_1}+k_{2,2}\varepsilon_4+\{l\}-\varepsilon_7}\mu_2\Pi.
\end{equation}

Let every rectangle $\Pi_k$ contain no more than $2^mQ^{[l]+\frac{\varepsilon_7}{2}}$ polynomials $P\in\Pl_m(Q,\mathbf{k}_1,\mathbf{k}_2,u)$. Then by inequalities (\ref{eq90}), (\ref{eq29}) and (\ref{eq92}) it follows that the measure of the set $L_{n+1,j}$ can be estimated as:
\[
\mu_2 L_{n+1,j} \ll Q^{-u-v_2+p_{2,1}+k_{2,2}\varepsilon_4-\frac{\varepsilon_7}{2}+[l]+\{l\}}\mu_2\Pi\ll Q^{-\frac{\varepsilon_7}{2}}\mu_2\Pi\leq \textstyle\frac{1}{12r}\mu_2\Pi,
\]
where $Q>Q_0$.

Now assume that there exists a rectangle $\Pi_k$ containing more than $2^mQ^{[l]+\frac{\varepsilon_7}{2}}$ points $(\alpha_1,\alpha_2)$, where $\alpha_1,\alpha_2$ are the roots of polynomial $P_j\in\Pl_m(Q,\mathbf{k}_1,\mathbf{k}_2,u)$. Using the calculations described in the previous case (see estimate (\ref{eq43})) and in Section \ref{sec_4} (see estimate (\ref{eq00050})) for every point $(x_1,x_2)\in\Pi_k$ we have:
\begin{equation}\label{eq93}
|P_j(x_1)|< Q^{-v_1+(m+1)(\varepsilon_4+\varepsilon_7)},\quad |P_j(x_2)|< Q^{u-p_{2,1}-k_{2,2}\varepsilon_4-\{l\}+(m+2)\varepsilon_4}.
\end{equation}

By Dirichlet's principle we can find at least $\left[Q^{\frac{\varepsilon_7}{2}}\right]+1$ from $\Pl_m(Q,\mathbf{k}_1,\mathbf{k}_2,u)$ contained in $\Pi_k$ such that their coefficients $a_m,\ldots, a_{m+1-[l]}$ coincide. Let us call them $P_1, \ldots, P_{\left[Q^{\frac{\varepsilon_7}{2}}\right]+1}$. Thus, let us consider the differences $R_{i,j}=P_i-P_j$, where $1\leq i<j\leq \left[Q^{\frac{\varepsilon_7}{2}}\right]+1$.

From the inequalities (\ref{eq93}), we obtain that at every point of the rectangle $\Pi_k$ the polynomials $R_{i,j}$ satisfy
\[
\begin{cases}
|R_{i,j}(x_1)|< 2Q^{-v_1+(m+1)(\varepsilon_4+\varepsilon_7)},\quad |R_{i,j}(x_2)|< 2Q^{u-p_{2,1}-k_{2,2}\varepsilon_4-\{l\}+(m+2)\varepsilon_4},\\
\deg R_{i,j} \leq m-[l]=m-v_2+p_{2,1}-u+k_{2,2}\varepsilon_4+\{l\}.
\end{cases}
\]

Assume that among polynomials $R_{i,j}$ we can find at least two polynomials without common roots and apply Lemma \ref{lm3} with $\tau_1= v_1-(m+1)(\varepsilon_4+\varepsilon_7)$, $\tau_2 = -u+p_{2,1}+k_{2,2}\cdot\varepsilon_4+\{l\}-(m+2)\varepsilon_4$, $\eta_1=\frac{u+v_1-p_{1,m_1}}{m_1}-\varepsilon_7$, $\eta_2=k_{2,2}\varepsilon_4+\{l\}$, so that we have
\[
\tau_1+1=v_1+1-(m+1)(\varepsilon_4+\varepsilon_7),\quad \tau_2+1=1-u+p_{2,1}+k_{2,2}\varepsilon_4+\{l\}-(m+2)\varepsilon_4,
\]
and repeating the arguments from the end of Section \ref{sec_4} we obtain
\[
2(\tau_1+1-\eta_1)\ge v_1+1-2(m+1)\varepsilon_4-2m\varepsilon_7,\quad 2(\tau_2+1-\eta_2)=2-2u+2p_{2,1}-2(m+2)\varepsilon_4.
\]
Substituting these expressions into (\ref{eq0}) for $\varepsilon_4=\textstyle\frac{1}{48(m+2)}$ and $\varepsilon_7=\textstyle\frac{1}{8(3m+1)}$ yields
\begin{multline*}
\tau_1+\tau_2+2+2(\tau_1+1-\eta_1)+2(\tau_2+1-\eta_2)\ge 2v_1+5+3p_{2,1}+k_{2,2}\varepsilon_4-3u+\{l\}
-\textstyle\frac14\ge\\ \ge 2n-2v_2+2p_{2,1}+2k_{2,2}\varepsilon_4-2u+\{l\}
+\textstyle\frac{7}{4}\ge 2(m-[l])-\{l\} +1+\textstyle\frac{3}{4}
\ge 2(m-[l])+\textstyle\textstyle\frac34.
\end{multline*}
This inequality contradicts to Lemma \ref{lm3} with $\delta =\frac34$. 

If among polynomials $R_{i,j}$, $1\leq i< j\leq \left[Q^{\frac{\varepsilon_7}{2}}\right]+1$ we can not find two polynomials without common roots then we use the arguments described in \cite{BernikGoetzeKukso14}.

This section concludes the proof of Lemma in case of irreducible polynomials. We have
\[
\mu_2 L_1 \leq \sum\limits_{2\leq i,j\leq n+1} \mu_2 L_{i,j}\leq (n-1)^2\cdot\textstyle\frac{1}{12r} \cdot\mu_2\Pi\leq \textstyle\frac{1}{12}\cdot \mu_2\Pi.
\]
Similarly we obtain $\mu_2 L_2\leq \textstyle\frac{1}{12}\cdot \mu_2\Pi$.

\subsection{The case of reducible polynomials}\label{sec_6}

Let us estimate the measure of the set $L_3$. Let a polynomial $P$ of degree $n$ be a product of several (not necessarily different) irreducible polynomials $P_1, P_2, \ldots, P_s$, $s\ge 2$, where $\deg P_i = n_i\ge 2$ and $n_1+\ldots+n_s=n$. Then by Lemma \ref{lm4} we have:
\[
H(P_1)\cdot H(P_2)\cdot\ldots\cdot H(P_s) \leq c_{9} H(P) \leq c_{9} Q.
\]
On the other hand, by the definition of height, we have $H(P_i)\ge 1$, and thus $H(P_i)\leq c_{9} Q=Q_1$, $i=1,\ldots,s$.

Denote by $L_{3}(k)$ a set of points $(x_1,x_2)\in\Pi$ such that there exists a polynomial $R\in\Pl_{k}(Q_1)$ satisfying the inequality:
\begin{equation}\label{eq73}
|R(x_1)R(x_2)|< h_n^2Q_1^{-k+\frac12}.
\end{equation}
If a polynomial $P\in\Pl_n(Q)$ satisfies the inequalities (\ref{eq5}) at a point $(x_1,x_2)\in\Pi$, we can write
\[
|P(x_1)P(x_2)|=|P_1(x_1)P_1(x_2)|\cdot\ldots\cdot|P_s(x_1)P_s(x_2)|\leq h_n^2Q^{-n+1}.
\]
Since $n= n_1 +\ldots+ n_s$ and $s\ge 2$, it is easy to see that at least one of the inequalities
\[
|P_i(x_1)P_i(x_2)| \leq h_n^2Q^{-n_i+\frac12},\quad i=1,\ldots, s,
\]
is satisfied at the point $(x_1,x_2)$. Hence, $(x_1,x_2)\in L_{3}(n_j)$ and we have 
\[
L_3\subset\bigcup\limits_{k=2}^{n-2}L_{3}(k).
\]

Let us estimate the measure of the set $L_{3}(k)$, $2\leq k\leq n-2$. Denote by $L_{3}^1(k,t)$ a set of points $(x_1,x_2)\in\Pi$ such that there exists a polynomial $P\in\Pl_{k}(Q_1)$ satisfying the inequalities:
\begin{equation}\label{eq75}
\begin{cases}
|P(x_1)|<h_{n}^2Q_1^t,\quad |P(x_2)|<h_{n}^2Q_1^{-k+1-t},\\
\min\limits_i\left\{|P'(\alpha_i)|\right\}<\delta_k Q_1,\quad x_i\in S(\alpha_i),i=1,2.
\end{cases}
\end{equation}
and by $L_{3}^2(k,t)$ a set of points $(x_1,x_2)\in\Pi$ such that there exists a polynomial $P\in\Pl_{k}(Q_1)$ satisfying the inequality:
\begin{equation}\label{eq76}
\begin{cases}
|P(x_1)|<h_{n}^2Q_1^t,\quad |P(x_2)|<h_{n}^2Q_1^{-k+\frac34-t},\\
|P'(\alpha_i)|>\delta_k Q_1,\quad x_i\in S(\alpha_i),\quad i=1,2.
\end{cases}
\end{equation}
By the definition of the set $L_3(k)$ it is easy to see that:
\[
L_3(k)\subset \left(\bigcup\limits_{i=0}^{2k}{L_3^1(k, 1-i/2)}\right)\cup\left(\bigcup\limits_{i=0}^{4k+1}{L_3^2(k, 1-i/4)}\right).
\]

The system (\ref{eq75}) is a system of the form (\ref{eq5}). Hence, as the polynomials $P\in \Pl_{k}(Q_1)$ are irreducible and $k<n$, we can apply the induction hypothesis to obtain the following estimate:
\begin{equation}\label{eq77}
\mu_2 L_3^1(k,t) < \textstyle\frac{1}{2^6n^2}\cdot \mu_2\Pi
\end{equation}
for $Q_1>Q_0$ and sufficiently small $\delta_k$.

Now let us estimate the measure of the set $L_3^2(k,t)$. From Lemma \ref{lm1} it follows that $L_{3}^2(k,t)$ is contained in a union $\bigcup\limits_{P\in\Pl_{k}(Q)}{\sigma_P(t)}$, where
\[
\sigma_{P}(t) := 
\left\{
(x_1,x_2)\in\Pi:
\begin{array}{ll}
|x_1-\alpha_1| \leq 2^{k-1}h_n^2\cdot Q^t\cdot|P'(\alpha_1)|^{-1},\\
|x_2-\alpha_2| \leq 2^{k-1}h_n^2\cdot Q^{-k+\frac34-t}\cdot|P'(\alpha_2)|^{-1}.
\end{array}
\right\}
\]
Let us estimate the value of the polynomial $P$ at a central point $\mathbf{d}$ of the square $\Pi$. A Taylor expansion of the polynomial $P$ can be written as follows:
\begin{equation}\label{eq78}
P(d_i)=P'(\alpha_i)(d_i-\alpha_i)+\textstyle\frac12 P''(\alpha_i)(d_i-\alpha_i)^2+\ldots +\textstyle\frac{1}{k!}\cdot P^{(k)}(\alpha_i)(d_i-\alpha_i)^k.
\end{equation}
If polynomial $P$ satisfy (\ref{eq76}) at point $(x_{0,1},x_{0,2})\in\Pi$ then:
\begin{equation}\label{eq79}
\begin{split}
|d_1-\alpha_1|&\leq |d_1-x_{0,1}|+|x_{0,1}-\alpha_1|\leq \mu_1 I_1 + 2^{k-1}h_n^2\delta_k^{-1}\cdot Q_1^{t-1},\\
|d_2-\alpha_2|&\leq |d_2-x_{0,2}|+|x_{0,2}-\alpha_2|\leq \mu_1 I_2 + 2^{k-1}h_n^2\delta_k^{-1}\cdot Q_1^{-k+\frac34-t-1}.
\end{split}
\end{equation}
Without loss of generality, let us assume that $t\ge -k+\frac34-t$. Then we can rewrite the estimates (\ref{eq79}) as follows:
\[
|d_1-\alpha_1|\leq 
\begin{cases}
c_{18}\cdot\mu_1 I_1,\quad t < 1-s_1,\\
c_{18}\cdot Q_1^{t-1},\quad 1-s_1 \leq t \leq 1,
\end{cases} \qquad |d_2-\alpha_2|\leq \mu_1 I_2.
\]
where $c_{18}=2^{k-1}h_n^2\delta_k^{-1}+c_{1,1}$. We mention that $\Pi = I_1\times I_2$, $\mu_1 I_i = c_{1,i}Q^{-s_i}$, $i=1,2$ and $s_1\leq s_2$.

Using these inequalities and expression (\ref{eq78}) allows us to write
\begin{equation}\label{eq80}
|P(d_1)|<
\begin{cases}
c_{19}\cdot  Q_1\cdot \mu_1 I_1,\quad t < 1-s_1,\\
c_{19}\cdot  Q_1^t,\quad 1-s_1 \leq t < 1,
\end{cases}\qquad  |P(d_2)|< c_{19} \cdot Q_1\cdot \mu_1 I_2.
\end{equation}

Fix a vector $\mathbf{A}_{1}=(a_k,\ldots,a_{2})$, where $a_k,\ldots,a_{2}$ will denote the coefficients of the polynomial $P\in \Pl_{k}(Q_1)$.
 Consider a subclass $\Pl_{k}(\mathbf{A}_{1})$ of polynomials $P$ which satisfy (\ref{eq76}) and have the same vector of coefficients $\mathbf{A}_{1}$. For $Q_1>Q_0$, the number of such classes can be estimated as follows
\begin{equation}\label{eq81}
\#\{\mathbf{A}_{1}\}=(2Q_1+1)^{k-1}< 2^{k} Q_1^{k-1}.
\end{equation}
Let us estimate the value $\#\Pl_{k}(\mathbf{A}_{1})$. Take a polynomial $P_0\in\Pl_{k}(\mathbf{A}_{1})$ and consider the difference between the polynomials $P_0$ and $P_j\in\Pl_{k}(\mathbf{A}_{1})$ at points $d_i$, $i=1,2$. By (\ref{eq80}), we have that:
\[
|P_0(d_1)-P_j(d_1)|=|(a_{0,1}-a_{j,1})d_1+(a_{0,0}-a_{j,0})|\leq \begin{cases}
2c_{19}\cdot  Q_1\mu_1 I_1,\quad t < 1-s_1,\\
2c_{19}\cdot  Q_1^t,\quad 1-s_1 \leq t \leq 1,
\end{cases}
\]
\[
|P_0(d_2)-P_j(d_2)|=|(a_{0,1}-a_{j,1})d_2+(a_{0,0}-a_{j,0})|\leq 2c_{19}\cdot  Q_1\mu_1 I_2.
\]
This implies that the number of different polynomials $P_j\in\Pl_{k}(\mathbf{A}_{1})$ does not exceed the number of integer solutions of the system
\begin{equation}\label{eq081}
|b_1d_i+b_0|\leq K_i,\quad i=1,2,
\end{equation}
where $K_2=2c_{19}\cdot  Q_1\mu_1 I_2$ and $K_1=
2c_{19}\cdot  Q_1\mu_1 I_1$ if $t < 1-s_1$ and $K_1=2c_{19}\cdot  Q_1^t$ if $1-s_1 \leq t \leq 1$. It is easy to see that $K_i\ge 2c_{19}\cdot Q_1^{1-s_1}>Q_1^{\varepsilon}>1$ for $Q_1>Q_0$. Thus, using the scheme described in Section \ref{sec_1} to solve the system (\ref{eq081}) we have
\[
\#\Pl_{k}(\mathbf{A}_{1})\leq
\begin{cases}
32\varepsilon^{-1}\cdot  Q_1^{2}\cdot \mu_2\Pi,\quad t < 1-s_1,\\
32\varepsilon^{-1}\cdot  Q_1^{t+1}\cdot \mu_1 I_2,\quad 1-s_1 \leq t \leq 1.
\end{cases}
\]

This estimate and the inequality (\ref{eq81}) mean that the number $N$ of polynomials $P\in \Pl_{k}(Q_1)$ satisfying the conditions (\ref{eq76}) can be estimated as follows:
\begin{equation}\label{eq82}
N\leq
\begin{cases}
2^{k+5}\varepsilon^{-1}\cdot  Q_1^{k+1}\cdot \mu_2\Pi,\quad t < 1-s_1,\\
2^{k+5}\varepsilon^{-1}\cdot  Q_1^{k+t}\cdot \mu_1 I_2,\quad 1-s_1 \leq t \leq 1.
\end{cases}
\end{equation}
On the other hand, the measure of the set $\sigma_{P}(t)$ satisfies the inequality
\begin{equation}\label{eq83}
\mu_2\sigma_{P}(t)\leq
\begin{cases}
2^{2k}h_n^4\delta_k^{-2}\cdot  Q_1^{-k-\frac54},\quad t < 1-s_1,\\
2^{2k}h_n^4\delta_k^{-2}\cdot  Q_1^{-k-\frac14-t}\cdot \mu_1 I_1,\quad 1-s_1 \leq t \leq 1.
\end{cases}
\end{equation}

Then, by estimates (\ref{eq82}) and (\ref{eq83}), for $Q_1>Q_0$ we can write
\begin{equation}\label{eq84}
\mu_2 L_{3}^2(k,t)\leq 2^{3k+5}\delta_k^{-2}h_n^4\varepsilon^{-1} Q_1^{-\frac14}\mu_2\Pi <\textstyle\frac{1}{2^7n^2}\cdot\mu_2\Pi.
\end{equation}

The inequalities (\ref{eq77}) and (\ref{eq84}) lead to the following estimate of the measure of the set $L_{3}(k)$:
\[
\mu_2 L_{3}(k) \leq \sum\limits_{i=0}^{2k}\mu_2L_3^1\left(k,1-i/2\right) + 
\sum\limits_{i=0}^{4k+1}\mu_2L_3^2\left(k,1-i/4\right)\leq \textstyle\frac{5+6k}{2^7n^2}\cdot\mu_2\Pi\leq \textstyle\frac{1}{12n}\cdot \mu_2\Pi.
\]
Therefore
\[
\mu_2 L_3\leq \sum\limits_{k=2}^{n-2}\mu_2 L_{3}(k)\leq \textstyle\frac{n-3}{12n}\cdot\mu_2\Pi\leq \textstyle\frac{1}{12}\cdot\mu_2\Pi.
\]
This proves Lemma \ref{lm5} in the case of reducible polynomials.

Thus, we have 
\[
\mu_2 L\leq \mu_2 L_1+\mu_2 L_2 +\mu_2 L_3 \leq \textstyle\frac{1}{4}\mu_2\Pi.
\]

\end{Proof}
\subsection{The final part of the proof}

The proof of Theorem 1 is going to be based on Lemma \ref{lm5}.
Consider a set
$B = \Pi \setminus L$. From Lemma \ref{lm5} it follows that
\begin{equation}\label{eq53}
\mu_2 B \ge \textstyle\frac34 \mu_2\Pi
\end{equation}
for $Q>Q_0$. 

It should be recalled that the value $h_n$ is defined in the beginning of the section 3 such that for every point $\mathbf{x}\in \Pi$ there exists a polynomial $P\in \Pl_n(Q)$ satisfying 
\[
|P(x_i)|\leq h_n Q^{-\frac{n-1}{2}},\quad i=1,2.
\]
Then, for every point $(x_{1,1},x_{1,2}) \in B$ there exists an irreducible polynomial
$P_1\in\mathcal{P}_n(Q)$ satisfying the system of inequalities
\[
\begin{cases}
|P_1(x_{1,i})|< h_nQ^{-\frac{n-1}{2}},\\
|P_1'(x_{1,i})|>\delta_n Q,\quad i=1,2.
\end{cases}
\]

Let $\alpha_i$, $x_{1,i}\in S(\alpha_i)$, $i=1,2$ be roots of the polynomial $P_1$. By Lemma \ref{lm1}, we have
\begin{equation}\label{eq55}
|x_{1,i}-\alpha_i|\leq nh_n\delta_n^{-1}Q^{-\frac{n+1}{2}},\quad i=1,2.
\end{equation}

We are going to choose a maximal system of points $\Gamma = (\boldsymbol{\gamma}_1,\ldots,\boldsymbol{\gamma}_t)$ satisfying the following conditions\\
1. $H(\boldsymbol{\gamma}_j)\leq Q, \deg(\boldsymbol{\gamma}_j)\leq n$;\\
2. Rectangles
\[
\sigma(\boldsymbol{\gamma}_j)=
\left\{
|x_i-\gamma_{j,i}|<nh_n\delta_n^{-1}Q^{-\frac{n+1}{2}},i=1,2\right\},\quad
j = \overline{1,t},
\]
do not intersect. 

Let us introduce an expanded rectangles
\begin{equation}\label{eq56}
\sigma_1(\boldsymbol{\gamma}_j)=
\left\{
|x_i-\gamma_{j,i}|<2nh_n\delta_n^{-1}Q^{-\frac{n+1}{2}},i=1,2\right\},\quad
j = \overline{1,t},
\end{equation}
and show that
\begin{equation}\label{eq57}
B\subset\bigcup_{j=1}^t \sigma_1(\boldsymbol{\gamma}_j).
\end{equation}

We obtain this by proving that for every point $(x_{1,1},x_{1,2}) \in B$ there exists a point $\boldsymbol{\gamma}_j\in\Gamma$ such that $(x_{1,1},x_{1,2})\in\sigma_1(\boldsymbol{\gamma}_j)$. 
Since $(x_{1,1},x_{1,2}) \in B$, there is a point $\boldsymbol{\alpha}=(\alpha_1,\alpha_2)$ such that the inequalities (\ref{eq55}) are true. Thus, either $\boldsymbol{\alpha}\in\Gamma$ and $(x_{1,1},x_{1,2})\in\sigma_1(\boldsymbol{\alpha})$, or there exists a point $\boldsymbol{\gamma}_j\in\Gamma$ satisfying
\[
|\alpha_i-\gamma_{j,i}|\leq nh_n\delta_n^{-1}Q^{-\frac{n+1}{2}},\quad i=1,2.
\]
Hence, $(x_{1,1},x_{1,2})\in\sigma_1(\boldsymbol{\gamma}_j)$.

In this case, by (\ref{eq53}),(\ref{eq56}) and (\ref{eq57}) we have:
\[
\textstyle\frac34 \mu_2\Pi \leq\mu_2 B\leq \sum\limits_{j=1}^t{\mu_2\sigma_1(\boldsymbol{\gamma}_j)}\leq t\cdot 2^4n^2h_n^2\delta_n^{-2}Q^{-n-1},
\]
which yields the estimate
\[
t \ge c_2 Q^{n+1}\mu_2\Pi.
\]

\section{Proof of Theorem 2}

The proof of Theorem 2 is based on the following Lemma.

\begin{lemma}\label{lm6}
For all {\it $\left(\textstyle\frac{v_1}{n-1},\textstyle\frac{v_2}{n-1}\right)$- ordinary} rectangles $\overline{\Pi}=I_1\times I_2$ such that:

1. $\mu_1 I_1=\mu_1 I_2=c_3Q^{-s}$, where $\textstyle\frac12 < s <\frac34$;

2. $\overline{\Pi}\cap\left\{(x_1,x_2)\in\R^2:|x_1-x_2|<\varepsilon\right\}=\emptyset$;

3. $c_3>c_0(n,\varepsilon,\mathbf{d})$, where $\mathbf{d}=(d_1,d_2)$ is the midpoint of $\overline{\Pi}$;

\noindent let $L=L(Q,\delta_n,\mathbf{v},\overline{\Pi})$ be the set of points $(x_1,x_2)\in\overline{\Pi}$ such that there exists a polynomial $P\in\Pl_n(Q)$ satisfying the following system of inequalities
\begin{equation}\label{eq61}
\begin{cases}
|P(x_i)|< h_n Q^{-v_i},\quad v_i>0,\\
\min\limits_i\left\{|P'(x_i)|\right\}<\delta_n Q,\\
v_1+v_2=n-1,\quad i=1,2.
\end{cases}
\end{equation}
Then for a sufficiently small constant $\delta_n<\delta_0(n,\varepsilon,\mathbf{d})$ and a sufficiently large $Q>Q_0(n,\varepsilon,\mathbf{v},s,\mathbf{d})$, the estimate
\[
\mu_2 L<\textstyle\frac14 \mu_2\overline{\Pi}
\]
holds.
\end{lemma}

\begin{Proof}
Lemma \ref{lm6} can be proved in the same way as Lemma \ref{lm5}, we only need to replace the base of induction. 

\begin{statement}
For all {\it $(\gamma_{2,1},\gamma_{2,2})$- ordinary} squares $\overline{\Pi}=I_1\times I_2$ under conditions 1 --- 3 let $L_{2,2}=L_{2,2}(Q,\delta_2,\boldsymbol{\gamma}_2,\overline{\Pi})$ be the set of points $(x_1,x_2)\in\overline{\Pi}$ such that there exists a polynomial $P\in\Pl_2(Q)$ satisfying the system of inequalities
\begin{equation}\label{eq62}
\begin{cases}
|P(x_i)|< h_2 Q^{-\gamma_{2,i}},\quad \gamma_{2,i}>0,\\
\min\limits_i\left\{|P'(x_i)|\right\}<\delta_2 Q,\quad i=1,2\\
\gamma_{2,i}+\gamma_{2,i}=1,\quad |b_2|>Q^{s-\frac12}.
\end{cases}
\end{equation}
Then for any $r>0$, $\delta_2\leq \delta_0(\varepsilon, r,\mathbf{d})$ and $Q>Q_0(n,\varepsilon,s,\boldsymbol{\gamma}_2,\mathbf{d})$, the estimate
\[
\mu_2 L_{2,2}<\textstyle\frac{1}{4r} \mu_2\overline{\Pi}
\]
holds.
\end{statement}

\begin{Proof}
Let $P$ be a polynomial of the form $P(t)=b_2t^2+b_1t+b_0$. Applying the same argument that we used to prove the Statement \ref{st1}, we obtain upper and lower bounds for the absolute value of the derivative $P'$ at roots $\alpha_1$, $\alpha_2$ and at points $x_1$, $x_2$, where $x_i\in S(\alpha_i)$, $i=1,2$:
\begin{equation}\label{eq62_1}
|P'(\alpha_i)|>\textstyle\frac34\cdot \varepsilon\cdot |b_2|,\quad
|P'(x_i)|\leq \left(|d_1|+|d_2|+1+\frac{\varepsilon}{4}\right)\cdot |b_2|.
\end{equation}
These estimates lead to the following inequality:
\[
|b_2|<4\delta_2\varepsilon^{-1}Q.
\]

From Lemma \ref{lm1} and the estimates (\ref{eq62}), (\ref{eq62_1}) it follows that the set $L_{2,2}$ is contained in a union $\bigcup\limits_{P\in\Pl_2(Q)}{\sigma_{P}}$, where
\begin{equation}\label{eq63}
\sigma_P=
\left\{(x_1,x_2)\in\overline{\Pi}:\quad
|x_i-\alpha_i|<2h_2\varepsilon^{-1}Q^{-\gamma_{2,i}}|b_2|^{-1},i=1,2
\right\}.
\end{equation}
Since the square $\overline{\Pi}$ is {\it $(\gamma_{2,1},\gamma_{2,2})$- ordinary} we have $|b_2|\ge Q^{s-\frac12}$ and
\[
\mu_2\sigma_P\leq 2^{4}h_2^2\varepsilon^{-2}Q^{-1}|b_2|^{-2}\leq c_3^2Q^{-2s}\leq\mu_2\overline{\Pi}
\]
with $c_3 > 4h_2\varepsilon^{-1}$ and $s>\textstyle\frac12$.

Then we can write the following estimate for the measure of the set $L_{2,2}$:
\[
\mu_2 L_{2,2}\leq\mu_2\bigcup\limits_{P\in\Pl_2(Q)}{\sigma_P}\leq\sum\limits_{P\in\Pl_2(Q)}{\mu_2\sigma_P}\leq 2^{4}h_2^2\varepsilon^{-2}Q^{-1}\sum\limits_{\substack{b_2,b_1,b_0\leq Q: \\ P(t)=b_2t^2+b_1t+b_0 \\ \sigma_{P}\neq\emptyset}}{|b_2|^{-2}}.
\]
As in the proof of Statement 1 of Lemma \ref{lm5}, we estimate the number of polynomials $P\in\Pl_2(Q)$ satisfying the system of inequalities (\ref{eq62}) at some point $(x_1,x_2)\in\overline{\Pi}$ for a fixed value of $b_2$.

Let us estimate the polynomial $P$ at the points $d_1,d_2$. From Taylor expansions and estimates (\ref{eq62_1}) we have
\begin{equation}\label{eq65}
|P(d_i)|\leq |P(x_i)|+c_{20}\cdot |b_2|\mu_1 I_i,
\end{equation}
for a sufficiently large $Q>Q_0$. Consider a system of equations
\begin{equation}\label{eq66}
\begin{cases}
b_2d_1^2+b_1d_1+b_0=l_1,\\
b_2d_2^2+b_1d_2+b_0=l_2,
\end{cases}
\end{equation}
in three variables $b_2,b_1,b_0\in\Z$, where $|l_i|\leq 2c_{20}\cdot\max\{1,|b_2|\mu_1 I_i\}$, $i=1,2$.

Let us estimate the number of possible pairs $(b_1,b_0)$ such that the system (\ref{eq66}) is satisfied for a fixed $b_2$. To obtain this estimate, we consider the system (\ref{eq66}) for two different combinations $b_2, b_{0,1},b_{0,0}$ and $b_2, b_{j,1},b_{j,0}$. Simple transformations lead to the following system of linear equations in two variables $b_{0,1}-b_{j,1}$ and $b_{0,0}-b_{j,0}$:
\begin{equation}\label{eq67}
\begin{cases}
(b_{0,1}-b_{j,1})d_1+(b_{0,0}-b_{j,0})=l_{0,1}-l_{j,1},\\
(b_{0,1}-b_{j,1})d_2+(b_{0,0}-b_{j,0})=l_{0,2}-l_{j,2}.
\end{cases}
\end{equation}
Since the determinant of this system does not vanish, we can use Cramer's rule to solve it. Using inequalities $|l_{0,i}-l_{j,i}|\leq 4c_{20}\cdot\max\{1,|b_2|\mu_1 I_i\}$ we estimate the determinants $\Delta_i$, $i=1,2$ as follows:
\[
|\Delta_i|\leq 8c_{20}\cdot \max\{1,|b_2|\mu_1 I_i\}.
\]
Thus
\[
|b_{0,i}-b_{j,i}|\leq\frac{|\Delta_i|}{|\Delta|}\leq 8c_{20}\varepsilon^{-1}\cdot\max\{1, |b_2|\mu_1 I_i\},
\]
and for a fixed $b_2$ the following estimate holds:
\begin{equation}\label{eq68}
\# (b_1,b_0)\leq
\begin{cases}
2^6c_{20}^2\varepsilon^{-2}|b_2|^2\mu_2\overline{\Pi},\quad|b_2|>c_3^{-1}Q^{s},\\
2^6c_{20}^2\varepsilon^{-2},\quad Q^{s-\frac12}<|b_2|<c_3^{-1}Q^{s}.
\end{cases}
\end{equation}

Depending on the absolute value $|b_2|$, let us consider the following two sets:
\[
L_{2,2}^{1}=\bigcup\limits_{\substack{P\in\Pl_2(Q),\\ c_3^{-1}Q^{s}<|b_2|<4\delta_2\varepsilon^{-1}Q}}\sigma_P,\qquad L_{2,2}^{2}=\bigcup\limits_{\substack{P\in\Pl_2(Q),\\ Q^{s-\frac12}<|b_2|<c_3^{-1}Q^{s}}}\sigma_P.
\]

{\it The set $L_{2,2}^1$}: In this case for $\delta_2 < 2^{-15}r^{-1}c_{20}^{-2}h_2^{-2}\varepsilon^5$ can be estimated as:
\[
\mu_2 L_{2,2}\leq 2^{10}c_{20}^2h_2^2\varepsilon^{-4} Q^{-1}\cdot 4\delta_2\varepsilon^{-1}Q\mu_2\overline{\Pi}<\textstyle\frac{1}{8r}\cdot\mu_2\overline{\Pi}.
\]

{\it The set $L_{2,2}^2$}: Consider the polynomials $P$ under condition $Q^{s-\frac12}<|b_2|<c_3^{-1}Q^{s}$. For every set $\sigma_P$ we define the expanded set:
\begin{equation}\label{eq069}
\sigma'_P=\left\{(x_1,x_2)\in\overline{\Pi}:\quad |x_i-\alpha_i|<2^5h_2\varepsilon^{-1}\sqrt{r}\cdot Q^{-\gamma_{2,i}}|b_2|^{-1},i=1,2\right\}.
\end{equation}
Let us prove that for $|b_2|<c_{21}\cdot Q^{\frac12}$, where $c_{21}=\varepsilon\left(2^5h_2\sqrt{r}\right)^{-1}\cdot \left(|d_1|+|d_2|+2\right)^{-1} $ this sets do not intersect.

Consider polynomials $P_j$, $j=1,2$ with roots $\alpha_{j,1},\alpha_{j,2}$ and leading coefficients $|b_{j,2}|<c_{21}\cdot Q^{\frac12}$. Without loss of generality we will assume $|b_{1,2}|<|b_{2,2}|$. Let there exists a point $(x_{0,1},x_{0,2})\in\sigma'_{P_1}\cap\sigma'_{P_2}$. Since $P_1$ and $P_2$ have no common roots, the resultant $R(P_1,P_2)$ doesn't vanish, and the following estimate holds:
\begin{equation}\label{eq69}
1\leq|R(P_1,P_2)|=|b_{1,2}|^2|b_{2,2}|^2|\alpha_{1,1}-\alpha_{2,1}||\alpha_{1,1}-\alpha_{2,2}||\alpha_{1,2}-\alpha_{2,1}||\alpha_{1,2}-\alpha_{2,2}|.
\end{equation}
By the estimates (\ref{eq069}) we have
\[
|\alpha_{1,i}-\alpha_{2,i}|\leq |\alpha_{1,i}-x_{0,i}|+|\alpha_{2,i}-x_{0,i}|<2^6h_2\varepsilon^{-1}\sqrt{r}\cdot Q^{-\gamma_{2,i}}|b_{1,2}|^{-1}.
\]
On the other hand for $Q>Q_0$ we get
\[
\begin{aligned}
|\alpha_{1,1}-\alpha_{2,2}|&\leq |\alpha_{1,1}|+|\alpha_{2,2}|\leq |d_1|+|d_2|+2,\\
|\alpha_{1,2}-\alpha_{2,1}|&\leq |\alpha_{1,2}|+|\alpha_{2,1}|\leq |d_1|+|d_2|+2.
\end{aligned}
\]
By substituting these inequalities to (\ref{eq69}) we obtain
\[
1\leq|R(P_1,P_2)|<2^{12}h_2^2\varepsilon^{-2}r\cdot \left(|d_1|+|d_2|+2\right)^2\cdot |b_{2,2}|^2\cdot Q^{-1}<\textstyle\frac14.
\]
This contradiction yields the following estimate
\[
\sum\limits_{\substack{P\in\Pl_2(Q),\\ Q^{s-\frac12}<|b_2|<c_{21}Q^{\frac12}}}\mu_2\sigma_P\leq \textstyle\frac{1}{16r}\cdot \sum\limits_{\substack{P\in\Pl_2(Q),\\ Q^{s-\frac12}<|b_2|<c_{21}Q^{\frac12}}}\mu_2\sigma'_P\leq \textstyle\frac{1}{16r}\cdot\mu_2\overline{\Pi}.
\]

Consider the case $|b_2|>c_{21}Q^{\frac12}$. Let $\Pl_{2}(Q,k)\subset\Pl_2(Q)$, $1\leq k\leq K=\left[\ln_2\left(\textstyle\frac{2-2s}{3-4s}\right)\right]+1$ be a subclass of polynomials defined as follows:
\[
\Pl_{2}(Q,k):=\left\{ P\in\Pl_2(Q): l_{k+1} \cdot Q^{\lambda_{k+1}}\leq |b_2|\leq l_k\cdot Q^{\lambda_k}\right\},
\]
where
\[ 
\begin{tabular}{lll}
$\lambda_1=s$, & $l_1= c_{3}^{-1}$, &\\
$\lambda_k = \lambda_{k-1} - (1-s)\cdot 2^{1-k}$, &  $l_k = \textstyle\frac{2^6c_{20}h_2\cdot \sqrt{rK\cdot l_{k-1}}}{\varepsilon^2 c_{3}}$ & for $2\leq k\leq K$,\\
$\lambda_{K+1} = \textstyle\frac12$, & $l_{K+1} =c_{21}$. &
\end{tabular}
\]
This equations give $\lambda_k = s - (1-s)\cdot\left(1-\textstyle\frac{1}{2^{k-1}}\right)$ for $2\leq k\leq K$.

Let us consider the following sets $L(k)=\bigcup\limits_{P\in \Pl_{2}(Q,k)}{\sigma_{P}}$ and estimate the measure of every one of them in the following way:
\[
\mu_2 L(k)=\sum\limits_{P\in \Pl_{2}(Q,k)}{\mu_2 \sigma_{P}}\leq \textstyle\frac{2^{10}h_2^2c_{20}^2}{\varepsilon^{4}}\cdot Q^{-1}\sum\limits_{l_{k+1} Q^{\lambda_{k+1}}\leq |b_2|\leq l_k Q^{\lambda_k}}{|b_2|^{-2}}
\leq \textstyle\frac{2^{10}h_2^2c_{20}^2l_k}{\varepsilon^{4}l_{k+1}^{2}}\cdot Q^{-1-2\lambda_{k+1}+\lambda_k}.
\]
Then for $k=1$ we obtain
\[
\mu_2 L(1)\leq \textstyle\frac{c_{3}^2}{16rK}\cdot Q^{-1-2s + 1-s+s}\leq \textstyle\frac{1}{16rK}\cdot c_{3}^2Q^{-2s}< \textstyle\frac{1}{16Kr}\cdot\mu_2\overline{\Pi};
\]
for $1\leq k\leq K-1$ we have
\[
\mu_2 L(k)\leq \textstyle\frac{c_{3}^2}{16rK}\cdot Q^{-1+s - (1-s)\cdot\left(1-\textstyle\frac{1}{2^{k-1}}\right)-2s + (1-s)\cdot\left(2-\textstyle\frac{1}{2^{k-1}}\right)}
\leq\textstyle\frac{1}{16rK}\cdot c_{3}^2 Q^{-2s}=\textstyle\frac{1}{16rK}\cdot\mu_2\overline{\Pi};
\]
and for $k=K$, $s<\frac34$ and $Q>Q_0$ we get
\[
\mu_2 L(K)\leq \textstyle\frac{2^6h_2^2\cdot l_K}{\varepsilon^{2}c_{20}^2}\cdot Q^{-2+s - (1-s)\cdot\left(1-\textstyle\frac{1}{2^{K-1}}\right)}\leq \textstyle\frac{2^4h_2^2\cdot l_K}{\varepsilon^{2}c_{21}^2}\cdot Q^{-3+2s+(1-s)\cdot\textstyle\frac{3-4s}{2-2s}}\leq\textstyle\frac{2^4h_2^2\cdot l_K}{\varepsilon^{2}c_{21}^2}\cdot Q^{-\frac32}< \textstyle\frac{1}{16rK}\cdot\mu_2\overline{\Pi}.
\]

Then, we obtain following estimate for the measure of the set $L_{2,2}^2$
\[
\mu_2 L_{2,2}^2\leq\sum\limits_{\substack{P\in\Pl_2(Q),\\  Q^{s-\frac12}<|b_2|<c_{21}Q^{\frac12}}}\mu_2\sigma_P+ \sum\limits_{1\leq k\leq K}\mu_2 L(k)\leq \textstyle\frac{1}{8r}\cdot\mu_2\overline{\Pi},
\]
and thus
\[
\mu_2 L_{2,2}\leq \mu_2 L_{2,2}^1 + \mu_2 L_{2,2}^2 \leq \textstyle\frac{1}{4r}\cdot\mu_2\overline{\Pi}.
\]

\end{Proof}

Now Lemma \ref{lm6} can be proved by repeating the proof of Lemma \ref{lm5}.
\end{Proof}

Theorem 2 is proved by applying the results of Lemma \ref{lm6} to the proof of Theorem 1.

\section{Proof of Theorem 3}

To prove Theorem 3 we are going to use the results of Theorem 1 and Theorem 2. For this purpose we need to consider the set $J\backslash D=\bigcup\limits_k J_k$, where $D:=\{x\in J:\quad |f(x)-x|<\textstyle\frac12\varepsilon\}$. It is easy to see, that for sufficiently small $\varepsilon$ the following estimate for the measure of the set $J\backslash D$ holds:
\[
\mu_1 (J\backslash D) \ge \textstyle\frac34\mu_1 J.
\]
Now for every strip $L_{J_k}(Q,\lambda)$ we have $L_{J_k}(Q,\lambda)\cap \left\{(x_1,x_2)\in\R^2:|x_1-x_2|<\varepsilon\right\}=\emptyset$.

Let us consider an interval $J_k=[a_k,b_k]$ and the strip $L_{J_k}(Q,\lambda)$ for a fixed $0<\lambda<\frac34$. Divide the strip $L_{J_k}(Q,\lambda)$ into segments 
\[
E_j:=\left\{(x_1,x_2)\in\R^2: x_1\in J_{k,j}, |x_2-f(x_1)|<\left(\textstyle\frac12+c_5\right)\cdot c_3Q^{-\lambda}\right\},
\]
where $J_{k,j}=[x_j, x_{j+1}]$, $x_j=x_{j-1}+c_3Q^{-\lambda}$, $x_0=a_k$ and $1\leq j\leq t_k$. The number of segments $E_j$ can be estimated as follows:
\[
t_k>\frac{\mu_1 J_k}{\mu_1 J_{k,j}}-1>\textstyle\frac12c_3^{-1}\mu_1 J_k\cdot Q^{\lambda}
\]
for $Q>Q_0$.

Let $\overline{f}_j=\textstyle\frac12\cdot\left(\max\limits_{x\in J_{k,j}}f(x)+\min\limits_{x\in J_{k,j}}f(x)\right)$. Consider the rectangles
\[
\Pi_j=\left\{(x_1,x_2)\in\R^2:x_1\in J_{k,j}, \left|x_2-\overline{f}_j\right|\leq \textstyle\frac12c_5c_3Q^{-\lambda}\right\}.
\]
By mean value theorem, since $f$ is continuous and differentiable function on every interval $J_{k,j}$ and $\sup\limits_{x\in J_{k,j}}{|f'(x)|}\leq\sup\limits_{x\in J}{|f'(x)|}:=c_5$, we obtain:
\[
\left|\max\limits_{x\in J_{k,j}}f(x)-\min\limits_{x\in J_{k,j}}f(x)\right|\leq |f'(\xi)|\cdot\mu_1 J_{k,j} < c_5c_3\cdot Q^{-\lambda}.
\]
It means that $\Pi_j\subset E_j$ for every $1\leq j\leq t_k$.

{\it Case 1}: $0< \lambda\leq\frac12$.

In this case, we apply the result of Theorem \ref{th1}.
From Theorem \ref{th1} it follows that every rectangle $\Pi_j$, $j=\overline{1,t_k}$, contains at least $c_{2}Q^{n+1-2\lambda}$ algebraic points of degree at most $n$ and height at most $Q$. Since we have $t_k>\textstyle\frac12c_3^{-1}\mu_1 J_k Q^{\lambda}$ and $\sum\limits_k \mu_1 J_k \ge \mu_1 J\backslash D >\textstyle\frac34\mu_1 J$, there must be at least $c_{6}Q^{n+1-\lambda}$ algebraic points $\mathbf{\alpha}\in L_J(Q,\lambda)\cap \A_n^2(Q)$.

{\it Case 2}: $\frac12< \lambda< \frac34$.

Theorem \ref{th2} will be used in that case. Let us count the number of {\it $\left(\frac12,\frac12\right)$- special} squares $\Pi_i$. By definition, a {\it $\left(\frac12,\frac12\right)$- special} square contains the points $(x_{0,1},x_{0,2})$ such that there exists a polynomial $P\in\mathcal{P}_2(Q)$ satisfying the system of inequalities
\begin{equation}\label{eq70}
\left\{
\begin{array}{lll}
|P(x_{0,i})|<h_2Q^{-\frac12},\quad i=1,2,\\
|b_2|\leq Q^{\lambda-\frac12}.
\end{array}
\right.
\end{equation}

Repeating the steps of the proof of Statement \ref{st1} from the beginning till inequality (\ref{eq12}), we obtain the following estimates:
\[
|P'(\alpha_1)|=|P'(\alpha_2)|>\textstyle\frac34\cdot \varepsilon\cdot |b_2|.
\]
Thus, by Lemma \ref{lm1} the set of points $(x_{1},x_{2})$ satisfying the system (\ref{eq70}) for a fixed polynomial $P$ is a subset of the following square:
\[
\sigma_P=\left\{(x_1,x_2)\in\R^2:|x_i-\alpha_i|\leq 2h_2\varepsilon^{-1}Q^{-\frac12}|b_2|^{-1},i=1,2\right\}.
\]
Let us estimate the number of squares $\Pi_j$, such that $\Pi_j\cap\sigma_P\neq \emptyset$. It is easy to see that the width of the strip $L_{J_k}(Q,\lambda)$ is smaller than the heights of the squares $\sigma_P$ for sufficiently large $c_3$. Hence, every square $\sigma_P$ intersects with at most $4h_2\varepsilon^{-1}c_3^{-1}Q^{\lambda-\frac12}|b_2|^{-1}$ squares $\Pi_j$. Therefore, the number $m_1$ of {\it $\left(\frac12,\frac12\right)$- special} squares $\Pi_i$ can be estimated as
\[
m_1\leq\sum\limits_{P\in\mathcal{P}_2(Q)}{4h_2\varepsilon^{-1}c_3^{-1}Q^{\lambda-\frac12}|b_2|^{-1}}\leq 4h_2\varepsilon^{-1}c_3^{-1}Q^{\lambda-\frac12}\sum\limits_{b_2,b_1,b_0}{|b_2|^{-1}}
\]
Now we need to estimate the number of polynomials $P\in\Pl_2(Q)$ satisfying the system of inequalities (\ref{eq70}) at some point $(x_1,x_2)\in L_{J_k}(Q,\lambda)$ for a fixed value of $b_2$. Since the function $f$ is continuously differentiable on the interval $J$, and $\sup\limits_{x\in J_k}{|f'(x)|}<c_{5}$, we get by the mean value theorem that
\[
\left|\max\limits_{x\in J_k}{f(x)}-\min\limits_{x\in J_k}{f(x)}\right|<c_{5}\cdot \mu_1 J_k,
\]
which implies that the set $L_{J_k}(Q,\lambda)$ is contained in a rectangle $\Pi=I_1\times I_2$, where $\mu_1 I_2 = c_5\mu_1 I_1=c_5\mu_1 J_k$.

Let us estimate the polynomial $P$ at the midpoint $(d_1,d_2)$ of the rectangle $\Pi$. Using the steps of the proof of Statement \ref{st1} we obtain
\[
|P(d_1)|\leq c_{22}\cdot |b_2|\mu_1 J_k,\quad |P(d_2)|\leq c_{22}c_5\cdot |b_2|\mu_1 J_k.
\]
and, hence, for a fixed value of $b_2$ the number of polynomials $P\in\Pl_2(Q)$ satisfying the system of inequalities (\ref{eq70}) at some point $(x_1,x_2)\in \Pi$ can be estimated as follows:
\[
\# (b_1,b_0)\leq 2^5c_{5}c_{22}^2\varepsilon^{-2}|b_2|^2\left(\mu_1 J_k\right)^2.
\]
Using this inequality we have:
\begin{equation}\label{eq71}
m_1\leq \textstyle\frac{2^7h_2c_{5}c_{22}^2\cdot \left(\mu_1 J_k\right)^2}{\varepsilon^{3}c_3}\cdot Q^{\lambda-\frac12}\sum\limits_{|b_2|< Q^{\lambda-\frac12}}{|b_2|}\leq \textstyle\frac{2^7h_2c_{5}c_{22}^2\cdot \left(\mu_1 J_k\right)^2}{\varepsilon^{3}c_3}\cdot Q^{3\lambda-\frac32}<\frac14c_3^{-1}\mu_1 J_{k}\cdot Q^{\lambda}<\textstyle\frac{t_k}{2}.
\end{equation}
for $\lambda<\frac34$ and $Q>Q_0$.
By (\ref{eq71}), it follows that the number of {\it $\left(\frac12,\frac12\right)$-ordinary} squares $\Pi_j$ doesn't exceed
\begin{equation}\label{eq72}
m_2\ge t_k-\textstyle\frac12t_k>\textstyle\frac12t_k.
\end{equation}

From Theorem 2 and the estimate (\ref{eq72}) it now follows that in the case $\frac12<\lambda< \frac34$, the strip $L_{J}(Q,\lambda)$ contains at least $c_{6}Q^{n+1-\lambda}$ algebraic points of degree at most $n$ and height at most $Q$.

{\bf Acknowledgements.} We would like to thank the referee for his thorough reports which
enabled us to remove a number of inaccuracies and improve the clarity of the paper.

This research was partly supported by SFB-701, Bielefeld University (Germany).

\end{document}